

\input amstex
\expandafter\ifx\csname mathdefs.tex\endcsname\relax
  \expandafter\gdef\csname mathdefs.tex\endcsname{}
\else \message{Hey!  Apparently you were trying to
  \string twice.   This does not make sense.} 
\errmessage{Please edit your file (probably \jobname.tex) and remove
any duplicate ``\string\input'' lines} \fi




\catcode`\X=12\catcode`\@=11

\def\n@wcount{\alloc@0\count\countdef\insc@unt}
\def\n@wwrite{\alloc@7\write\chardef\sixt@@n}
\def\n@wread{\alloc@6\read\chardef\sixt@@n}
\def\r@s@t{\relax}\def\v@idline{\par}\def\@mputate#1/{#1}
\def\l@c@l#1X{\firstpart.#1}\def\gl@b@l#1X{#1}\def\t@d@l#1X{{}}

\def\crossrefs#1{\ifx\all#1\let\tr@ce=\all\else\def\tr@ce{#1,}\fi
   \n@wwrite\cit@tionsout\openout\cit@tionsout=\jobname.cit 
   \write\cit@tionsout{\tr@ce}\expandafter\setfl@gs\tr@ce,}
\def\setfl@gs#1,{\def\@{#1}\ifx\@\empty\let\next=\relax
   \else\let\next=\setfl@gs\expandafter\xdef
   \csname#1tr@cetrue\endcsname{}\fi\next}
\def\m@ketag#1#2{\expandafter\n@wcount\csname#2tagno\endcsname
     \csname#2tagno\endcsname=0\let\tail=\all\xdef\all{\tail#2,}
   \ifx#1\l@c@l\let\tail=\r@s@t\xdef\r@s@t{\csname#2tagno\endcsname=0\tail}\fi
   \expandafter\gdef\csname#2cite\endcsname##1{\expandafter
     \ifx\csname#2tag##1\endcsname\relax?\else\csname#2tag##1\endcsname\fi
     \expandafter\ifx\csname#2tr@cetrue\endcsname\relax\else
     \write\cit@tionsout{#2tag ##1 cited on page \folio.}\fi}
   \expandafter\gdef\csname#2page\endcsname##1{\expandafter
     \ifx\csname#2page##1\endcsname\relax?\else\csname#2page##1\endcsname\fi
     \expandafter\ifx\csname#2tr@cetrue\endcsname\relax\else
     \write\cit@tionsout{#2tag ##1 cited on page \folio.}\fi}
   \expandafter\gdef\csname#2tag\endcsname##1{\expandafter
      \ifx\csname#2check##1\endcsname\relax
      \expandafter\xdef\csname#2check##1\endcsname{}%
      \else\immediate\write16{Warning: #2tag ##1 used more than once.}\fi
      \multit@g{#1}{#2}##1/X%
      \write\t@gsout{#2tag ##1 assigned number \csname#2tag##1\endcsname\space
      on page \number\count0.}%
   \csname#2tag##1\endcsname}}
\def\multit@g#1#2#3/#4X{\def\t@mp{#4}\ifx\t@mp\empty%
      \global\advance\csname#2tagno\endcsname by 1 
      \expandafter\xdef\csname#2tag#3\endcsname
      {#1\number\csname#2tagno\endcsnameX}%
   \else\expandafter\ifx\csname#2last#3\endcsname\relax
      \expandafter\n@wcount\csname#2last#3\endcsname
      \global\advance\csname#2tagno\endcsname by 1 
      \expandafter\xdef\csname#2tag#3\endcsname
      {#1\number\csname#2tagno\endcsnameX}
      \write\t@gsout{#2tag #3 assigned number \csname#2tag#3\endcsname\space
      on page \number\count0.}\fi
   \global\advance\csname#2last#3\endcsname by 1
   \def\t@mp{\expandafter\xdef\csname#2tag#3/}%
   \expandafter\t@mp\@mputate#4\endcsname
   {\csname#2tag#3\endcsname\lastpart{\csname#2last#3\endcsname}}\fi}
\def\t@gs#1{\def\all{}\m@ketag#1e\m@ketag#1s\m@ketag\t@d@l p
   \m@ketag\gl@b@l r \n@wread\t@gsin
   \openin\t@gsin=\jobname.tgs \re@der \closein\t@gsin
   \n@wwrite\t@gsout\openout\t@gsout=\jobname.tgs }
\outer\def\localtags{\t@gs\l@c@l}
\outer\def\globaltags{\t@gs\gl@b@l}
\outer\def\newlocaltag#1{\m@ketag\l@c@l{#1}}
\outer\def\newglobaltag#1{\m@ketag\gl@b@l{#1}}

\newif\ifpr@ 
\def\m@kecs #1tag #2 assigned number #3 on page #4.%
   {\expandafter\gdef\csname#1tag#2\endcsname{#3}
   \expandafter\gdef\csname#1page#2\endcsname{#4}
   \ifpr@\expandafter\xdef\csname#1check#2\endcsname{}\fi}
\def\re@der{\ifeof\t@gsin\let\next=\relax\else
   \read\t@gsin to\t@gline\ifx\t@gline\v@idline\else
   \expandafter\m@kecs \t@gline\fi\let \next=\re@der\fi\next}
\def\pretags#1{\pr@true\pret@gs#1,,}
\def\pret@gs#1,{\def\@{#1}\ifx\@\empty\let\n@xtfile=\relax
   \else\let\n@xtfile=\pret@gs \openin\t@gsin=#1.tgs \message{#1} \re@der 
   \closein\t@gsin\fi \n@xtfile}

\newcount\sectno\sectno=0\newcount\subsectno\subsectno=0
\newif\ifultr@local \def\ultralocal{\ultr@localtrue}
\def\firstpart{\number\sectno}
\def\lastpart#1{\ifcase#1 \or a\or b\or c\or d\or e\or f\or g\or h\or 
   i\or k\or l\or m\or n\or o\or p\or q\or r\or s\or t\or u\or v\or w\or 
   x\or y\or z \fi}

\def\resetall{\global\advance\sectno by 1\subsectno=0
   \gdef\firstpart{\number\sectno}\r@s@t}
\def\resetsub{\global\advance\subsectno by 1
   \gdef\firstpart{\number\sectno.\number\subsectno}\r@s@t}
\def\newsection#1\par{\resetall\vskip0pt plus.3\vsize\penalty-250
   \vskip0pt plus-.3\vsize\bigskip\bigskip
   \message{#1}\leftline{\bf#1}\nobreak\bigskip}
\def\subsection#1\par{\ifultr@local\resetsub\fi
   \vskip0pt plus.2\vsize\penalty-250\vskip0pt plus-.2\vsize
   \bigskip\smallskip\message{#1}\leftline{\bf#1}\nobreak\medskip}

\def\t@gsoff#1,{\def\@{#1}\ifx\@\empty\let\next=\relax\else\let\next=\t@gsoff
   \def\@@{p}\ifx\@\@@\else
   \expandafter\gdef\csname#1cite\endcsname##1{\zeigen{##1}}
   \expandafter\gdef\csname#1page\endcsname##1{?}
   \expandafter\gdef\csname#1tag\endcsname##1{\zeigen{##1}}\fi\fi\next}
\def\verbatimtags{\ifx\all\relax\else\expandafter\t@gsoff\all,\fi}
\def\zeigen#1{\hbox{$\langle$}#1\hbox{$\rangle$}}

\def\(#1){\edef\dot@g{\ifmmode\ifinner(\hbox{\noexpand\etag{#1}})
   \else\noexpand\eqno(\hbox{\noexpand\etag{#1}})\fi
   \else(\noexpand\ecite{#1})\fi}\dot@g}

\newif\ifbr@ck
\def\eat#1{}
\def\[#1]{\br@cktrue[\br@cket#1'X]}
\def\br@cket#1'#2X{\def\temp{#2}\ifx\temp\empty\let\next\eat
   \else\let\next\br@cket\fi
   \ifbr@ck\br@ckfalse\br@ck@t#1,X\else\br@cktrue#1\fi\next#2X}
\def\br@ck@t#1,#2X{\def\temp{#2}\ifx\temp\empty\let\neext\eat
   \else\let\neext\br@ck@t\def\temp{,}\fi
   \def\teemp{#1}\ifx\teemp\empty\else\rcite{#1}\fi\temp\neext#2X}
\def\resetbr@cket{\gdef\[##1]{[\rtag{##1}]}}
\def\references{\resetbr@cket\newsection References\par}

\newtoks\symb@ls\newtoks\s@mb@ls\newtoks\p@gelist\n@wcount\ftn@mber
    \ftn@mber=1\newif\ifftn@mbers\ftn@mbersfalse\newif\ifbyp@ge\byp@gefalse
\def\defm@rk{\ifftn@mbers\n@mberm@rk\else\symb@lm@rk\fi}
\def\n@mberm@rk{\xdef\m@rk{{\the\ftn@mber}}%
    \global\advance\ftn@mber by 1 }
\def\rot@te#1{\let\temp=#1\global#1=\expandafter\r@t@te\the\temp,X}
\def\r@t@te#1,#2X{{#2#1}\xdef\m@rk{{#1}}}
\def\b@@st#1{{$^{#1}$}}\def\str@p#1{#1}
\def\symb@lm@rk{\ifbyp@ge\rot@te\p@gelist\ifnum\expandafter\str@p\m@rk=1 
    \s@mb@ls=\symb@ls\fi\write\f@nsout{\number\count0}\fi \rot@te\s@mb@ls}
\def\byp@ge{\byp@getrue\n@wwrite\f@nsin\openin\f@nsin=\jobname.fns 
    \n@wcount\currentp@ge\currentp@ge=0\p@gelist={0}
    \re@dfns\closein\f@nsin\rot@te\p@gelist
    \n@wread\f@nsout\openout\f@nsout=\jobname.fns }
\def\m@kelist#1X#2{{#1,#2}}
\def\re@dfns{\ifeof\f@nsin\let\next=\relax\else\read\f@nsin to \f@nline
    \ifx\f@nline\v@idline\else\let\t@mplist=\p@gelist
    \ifnum\currentp@ge=\f@nline
    \global\p@gelist=\expandafter\m@kelist\the\t@mplistX0
    \else\currentp@ge=\f@nline
    \global\p@gelist=\expandafter\m@kelist\the\t@mplistX1\fi\fi
    \let\next=\re@dfns\fi\next}
\def\symbols#1{\symb@ls={#1}\s@mb@ls=\symb@ls} 
\def\bigsymbol{\textstyle}
\symbols{\bigsymbol\ast,\dagger,\ddagger,\sharp,\flat,\natural,\star}
\def\ftnumbers{\ftn@mberstrue} \def\ftsymbols{\ftn@mbersfalse}
\def\paginal{\byp@ge} \def\resetftnumbers{\ftn@mber=1}
\def\ftnote#1{\defm@rk\expandafter\expandafter\expandafter\footnote
    \expandafter\b@@st\m@rk{#1}}

\long\def\jump#1\endjump{}
\def\ssum{\mathop{\lower .1em\hbox{$\textstyle\Sigma$}}\nolimits}

\def\qed{\nobreak\kern 1em \vrule height .5em width .5em depth 0em}
\def\newneq{\hbox{\rlap{\hbox to 1\wd9{\hss$=$\hss}}\raise .1em 
   \hbox to 1\wd9{\hss$\scriptscriptstyle/$\hss}}}
\def\subsetne{\setbox9 = \hbox{$\subset$}\mathrel{\hbox{\rlap
   {\lower .4em \newneq}\raise .13em \hbox{$\subset$}}}}
\def\supsetne{\setbox9 = \hbox{$\subset$}\mathrel{\hbox{\rlap
   {\lower .4em \newneq}\raise .13em \hbox{$\supset$}}}}

\def\vbar{\mathchoice{\vrule height6.3ptdepth-.5ptwidth.8pt\kern-.8pt}
   {\vrule height6.3ptdepth-.5ptwidth.8pt\kern-.8pt}
   {\vrule height4.1ptdepth-.35ptwidth.6pt\kern-.6pt}
   {\vrule height3.1ptdepth-.25ptwidth.5pt\kern-.5pt}}
\def\f@dge{\mathchoice{}{}{\mkern.5mu}{\mkern.8mu}}
\def\b@c#1#2{{\rm \mkern#2mu\vbar\mkern-#2mu#1}}
\def\b@b#1{{\rm I\mkern-3.5mu #1}}
\def\b@a#1#2{{\rm #1\mkern-#2mu\f@dge #1}}
\def\bb#1{{\count4=`#1 \advance\count4by-64 \ifcase\count4\or\b@a A{11.5}\or
   \b@b B\or\b@c C{5}\or\b@b D\or\b@b E\or\b@b F \or\b@c G{5}\or\b@b H\or
   \b@b I\or\b@c J{3}\or\b@b K\or\b@b L \or\b@b M\or\b@b N\or\b@c O{5} \or
   \b@b P\or\b@c Q{5}\or\b@b R\or\b@a S{8}\or\b@a T{10.5}\or\b@c U{5}\or
   \b@a V{12}\or\b@a W{16.5}\or\b@a X{11}\or\b@a Y{11.7}\or\b@a Z{7.5}\fi}}

\catcode`\X=11 \catcode`\@=12

\expandafter\ifx\csname citeadd.tex\endcsname\relax
\expandafter\gdef\csname citeadd.tex\endcsname{}
\else \message{Hey!  Apparently you were trying to
\string twice.   This does not make sense.} 
\errmessage{Please edit your file (probably \jobname.tex) and remove
any duplicate ``\string\input'' lines} \fi

\sectno=-1   
\localtags
\ifx\shlhetal\undefinedcontrolsequence\let\shlhetal\relax\fi
\NoBlackBoxes
\define\mr{\medskip\roster}
\define\nl{\newline}
\define\sn{\smallskip\noindent}
\define\mn{\medskip\noindent}
\define\bn{\bigskip\noindent}
\define\ub{\underbar}

\define\ermn{\endroster\medskip\noindent}
\define\dbca{\dsize\bigcap}
\define\dbcu{\dsize\bigcup}
\documentstyle {amsppt}
\topmatter
\title{{Non-Existence of Universal Members in Classes of Abelian Groups}\\
 Sh622} \endtitle
\rightheadtext{Non-Existence of Universal Members}
\author {Saharon Shelah \thanks {\null\newline This research was partially
supported by the German-Israel Foundation for Science \null\newline
I would like to thank Alice Leonhardt for the beautiful typing. \null\newline
First Typed - 96/Oct/18 \null\newline
Latest Revision - 98/Aug/26} \endthanks} \endauthor
\affil{Institute of Mathematics\\
 The Hebrew University\\
 Jerusalem, Israel
 \medskip
 Rutgers University\\
 Mathematics Department\\
 New Brunswick, NJ  USA
 \medskip
 The Fields Institute\\
 Toronto, Ontario\\
 Canada} \endaffil
\abstract {We prove that if $\mu^+ < \lambda = \text{ cf}(\lambda) <
\mu^{\aleph_0}$, then there is no universal reduced torsion free abelian
group.  Similarly if $\aleph_0 < \lambda < 2^{\aleph_0}$.  We also prove
that if $2^{\aleph_0} < \mu^+ < \lambda = \text{ cf}(\lambda) < 
\mu^{\aleph_0}$, then there is no universal reduced separable abelian
$p$-group in $\lambda$.  (Note: both results fail if (a) $\lambda =
\lambda^{\aleph_0}$ or if (b) $\lambda$ is strong limit, cf$(\mu) = \aleph_0 <
\mu$).} \endabstract
\endtopmatter
\document  

\expandafter\ifx\csname alice2jlem.tex\endcsname\relax
  \expandafter\gdef\csname alice2jlem.tex\endcsname{}
\else \message{Hey!  Apparently you were trying to
\string  twice.   This does not make sense.}
\errmessage{Please edit your file (probably \jobname.tex) and remove
any duplicate ``\string\input'' lines} \fi

\expandafter\ifx\csname bib4plain.tex\endcsname\relax
  \expandafter\gdef\csname bib4plain.tex\endcsname{}
\else \message{Hey!  Apparently you were trying to \string twice.   This does not make sense.}
\errmessage{Please edit your file (probably \jobname.tex) and remove
any duplicate ``\string\input'' lines} \fi

\def\renewcommand{\newcommand}	       
\edef\cite{\the\catcode`@}%
\catcode`@ = 11
\let\@oldatcatcode = \cite
\chardef\@letter = 11
\chardef\@other = 12
%
%
%
%
\def\@innerdef#1#2{\edef#1{\expandafter\noexpand\csname #2\endcsname}}%
%
%
\@innerdef\@innernewcount{newcount}%
\@innerdef\@innernewdimen{newdimen}%
\@innerdef\@innernewif{newif}%
\@innerdef\@innernewwrite{newwrite}%
%
%
%
\def\@gobble#1{}%
%
%
%
\ifx\inputlineno\@undefined
   \let\@linenumber = \empty 
\else
   \def\@linenumber{\the\inputlineno:\space}%
\fi
%
%
%
\def\@futurenonspacelet#1{\def\cs{#1}%
   \afterassignment\@stepone\let\@nexttoken=
}%
\begingroup 
\def\\{\global\let\@stoken= }%
\\ 
\endgroup
\def\@stepone{\expandafter\futurelet\cs\@steptwo}%
\def\@steptwo{\expandafter\ifx\cs\@stoken\let\@@next=\@stepthree
   \else\let\@@next=\@nexttoken\fi \@@next}%
\def\@stepthree{\afterassignment\@stepone\let\@@next= }%
%
%
%
\def\@getoptionalarg#1{%
   \let\@optionaltemp = #1%
   \let\@optionalnext = \relax
   \@futurenonspacelet\@optionalnext\@bracketcheck
}%
%
%
\def\@bracketcheck{%
   \ifx [\@optionalnext
      \expandafter\@@getoptionalarg
   \else
      \let\@optionalarg = \empty
      \expandafter\@optionaltemp
   \fi
}%
\def\@@getoptionalarg[#1]{%
   \def\@optionalarg{#1}%
   \@optionaltemp
}%
%
%
%
\def\@nnil{\@nil}%
\def\@fornoop#1\@@#2#3{}%
\def\@for#1:=#2\do#3{%
   \edef\@fortmp{#2}%
   \ifx\@fortmp\empty \else
      \expandafter\@forloop#2,\@nil,\@nil\@@#1{#3}%
   \fi
}%
\def\@forloop#1,#2,#3\@@#4#5{\def#4{#1}\ifx #4\@nnil \else
       #5\def#4{#2}\ifx #4\@nnil \else#5\@iforloop #3\@@#4{#5}\fi\fi
}%
\def\@iforloop#1,#2\@@#3#4{\def#3{#1}\ifx #3\@nnil
       \let\@nextwhile=\@fornoop \else
      #4\relax\let\@nextwhile=\@iforloop\fi\@nextwhile#2\@@#3{#4}%
}%
%
%
%
\@innernewif\if@fileexists
\def\@testfileexistence{\@getoptionalarg\@finishtestfileexistence}%
\def\@finishtestfileexistence#1{%
   \begingroup
      \def\extension{#1}%
      \immediate\openin0 =
         \ifx\@optionalarg\empty\jobname\else\@optionalarg\fi
         \ifx\extension\empty \else .#1\fi
         \space
      \ifeof 0
         \global\@fileexistsfalse
      \else
         \global\@fileexiststrue
      \fi
      \immediate\closein0
   \endgroup
}%
%
%
%
%
\def\bibliographystyle#1{%
   \@readauxfile
   \@writeaux{\string\bibstyle{#1}}%
}%
\let\bibstyle = \@gobble
%
%
\let\bblfilebasename = \jobname
\def\bibliography#1{%
   \@readauxfile
   \@writeaux{\string\bibdata{#1}}%
   \@testfileexistence[\bblfilebasename]{bbl}%
   \if@fileexists
      \nobreak
      \@readbblfile
   \fi
}%
\let\bibdata = \@gobble
%
%
\def\nocite#1{%
   \@readauxfile
   \@writeaux{\string\citation{#1}}%
}%
\@innernewif\if@notfirstcitation
%
%
\def\cite{\@getoptionalarg\@cite}%
%
%
\def\@cite#1{%
   \let\@citenotetext = \@optionalarg
   \printcitestart
   \nocite{#1}%
   \@notfirstcitationfalse
   \@for \@citation :=#1\do
   {%
      \expandafter\@onecitation\@citation\@@
   }%
   \ifx\empty\@citenotetext\else
      \printcitenote{\@citenotetext}%
   \fi
   \printcitefinish
}%
\def\@onecitation#1\@@{%
   \if@notfirstcitation
      \printbetweencitations
   \fi
   \expandafter \ifx \csname\@citelabel{#1}\endcsname \relax
      \if@citewarning
         \message{\@linenumber Undefined citation `#1'.}%
      \fi
      \expandafter\gdef\csname\@citelabel{#1}\endcsname{%
\strut
\vadjust{\vskip-\dp\strutbox
\vbox to 0pt{\vss\parindent0cm \leftskip=\hsize 
\advance\leftskip3mm
\advance\hsize 4cm\strut\openup-4pt 
\rightskip 0cm plus 1cm minus 0.5cm ?  #1 ?\strut}}
         {\tt
            \escapechar = -1
            \nobreak\hskip0pt
            \expandafter\string\csname#1\endcsname
            \nobreak\hskip0pt
         }%
      }%
   \fi
   \csname\@citelabel{#1}\endcsname
   \@notfirstcitationtrue
}%
%
%
\def\@citelabel#1{b@#1}%
%
%
\def\@citedef#1#2{\expandafter\gdef\csname\@citelabel{#1}\endcsname{#2}}%
%
%
%
\def\@readbblfile{%
   \ifx\@itemnum\@undefined
      \@innernewcount\@itemnum
   \fi
   \begingroup
      \def\begin##1##2{%
         \setbox0 = \hbox{\biblabelcontents{##2}}%
         \biblabelwidth = \wd0
      }%
      \def\end##1{}
      %
      %
      \@itemnum = 0
      \def\bibitem{\@getoptionalarg\@bibitem}%
      \def\@bibitem{%
         \ifx\@optionalarg\empty
            \expandafter\@numberedbibitem
         \else
            \expandafter\@alphabibitem
         \fi
      }%
      \def\@alphabibitem##1{%
         \expandafter \xdef\csname\@citelabel{##1}\endcsname {\@optionalarg}%
         \ifx\biblabelprecontents\@undefined
            \let\biblabelprecontents = \relax
         \fi
         \ifx\biblabelpostcontents\@undefined
            \let\biblabelpostcontents = \hss
         \fi
         \@finishbibitem{##1}%
      }%
      \def\@numberedbibitem##1{%
         \advance\@itemnum by 1
         \expandafter \xdef\csname\@citelabel{##1}\endcsname{\number\@itemnum}%
         \ifx\biblabelprecontents\@undefined
            \let\biblabelprecontents = \hss
         \fi
         \ifx\biblabelpostcontents\@undefined
            \let\biblabelpostcontents = \relax
         \fi
         \@finishbibitem{##1}%
      }%
      \def\@finishbibitem##1{%
         \biblabelprint{\csname\@citelabel{##1}\endcsname}%
         \@writeaux{\string\@citedef{##1}{\csname\@citelabel{##1}\endcsname}}%
         \ignorespaces
      }%
      %
      %
      \let\em = \bblem
      \let\newblock = \bblnewblock
      \let\sc = \bblsc
      \frenchspacing
      \clubpenalty = 4000 \widowpenalty = 4000
      \tolerance = 10000 \hfuzz = .5pt
      \everypar = {\hangindent = \biblabelwidth
                      \advance\hangindent by \biblabelextraspace}%
      \bblrm
      \parskip = 1.5ex plus .5ex minus .5ex
      \biblabelextraspace = .5em
      \bblhook
      \input \bblfilebasename.bbl
   \endgroup
}%
%
%
\@innernewdimen\biblabelwidth
\@innernewdimen\biblabelextraspace
%
%
%
\def\biblabelprint#1{%
   \noindent
   \hbox to \biblabelwidth{%
      \biblabelprecontents
      \biblabelcontents{#1}%
      \biblabelpostcontents
   }%
   \kern\biblabelextraspace
}%
%
%
%
\def\biblabelcontents#1{{\bblrm [#1]}}%
%
%
\def\bblrm{\rm}%
%
%
\def\bblem{\it}%
%
%
\def\bblsc{\ifx\@scfont\@undefined
              \font\@scfont = cmcsc10
           \fi
           \@scfont
}%
%
%
\def\bblnewblock{\hskip .11em plus .33em minus .07em }%
%
%
\let\bblhook = \empty
%
%
%
\def\printcitestart{[}
\def\printcitefinish{]}
\def\printbetweencitations{, }
\def\printcitenote#1{, #1}
%
%
%
\let\citation = \@gobble
%
%
%
\@innernewcount\@numparams
%
%
\def\newcommand#1{%
   \def\@commandname{#1}%
   \@getoptionalarg\@continuenewcommand
}%
%
%
\def\@continuenewcommand{%
   \@numparams = \ifx\@optionalarg\empty 0\else\@optionalarg \fi \relax
   \@newcommand
}%
%
%
\def\@newcommand#1{%
   \def\@startdef{\expandafter\edef\@commandname}%
   \ifnum\@numparams=0
      \let\@paramdef = \empty
   \else
      \ifnum\@numparams>9
         \errmessage{\the\@numparams\space is too many parameters}%
      \else
         \ifnum\@numparams<0
            \errmessage{\the\@numparams\space is too few parameters}%
         \else
            \edef\@paramdef{%
               \ifcase\@numparams
                  \empty  No arguments.
               \or ####1%
               \or ####1####2%
               \or ####1####2####3%
               \or ####1####2####3####4%
               \or ####1####2####3####4####5%
               \or ####1####2####3####4####5####6%
               \or ####1####2####3####4####5####6####7%
               \or ####1####2####3####4####5####6####7####8%
               \or ####1####2####3####4####5####6####7####8####9%
               \fi
            }%
         \fi
      \fi
   \fi
   \expandafter\@startdef\@paramdef{#1}%
}%
%
%
%
%
\def\@readauxfile{%
   \if@auxfiledone \else 
      \global\@auxfiledonetrue
      \@testfileexistence{aux}%
      \if@fileexists
         \begingroup
            \endlinechar = -1
            \catcode`@ = 11
            \input \jobname.aux
         \endgroup
      \else
         \message{\@undefinedmessage}%
         \global\@citewarningfalse
      \fi
      \immediate\openout\@auxfile = \jobname.aux
   \fi
}%
%
%
\newif\if@auxfiledone
\ifx\noauxfile\@undefined \else \@auxfiledonetrue\fi
%
%
%
%
\@innernewwrite\@auxfile
\def\@writeaux#1{\ifx\noauxfile\@undefined \write\@auxfile{#1}\fi}%
%
%
%
\ifx\@undefinedmessage\@undefined
   \def\@undefinedmessage{No .aux file; I won't give you warnings about
                          undefined citations.}%
\fi
%
%
\@innernewif\if@citewarning
\ifx\noauxfile\@undefined \@citewarningtrue\fi
%
%
%
\catcode`@ = \@oldatcatcode


\def\widestnumber#1#2{}

\def\rm{\fam0 \tenrm}

\def\fakesubhead#1\endsubhead{\bigskip\noindent{\bf#1}\par}


%
%
%

%

\font\textrsfs=rsfs10
\font\scriptrsfs=rsfs7
\font\scriptscriptrsfs=rsfs5

\newfam\rsfsfam
\textfont\rsfsfam=\textrsfs
\scriptfont\rsfsfam=\scriptrsfs
\scriptscriptfont\rsfsfam=\scriptscriptrsfs

\edef\oldcatcodeofat{\the\catcode`\@}
\catcode`\@11

\def\Cal@@#1{\noaccents@ \fam \rsfsfam #1}

\catcode`\@\oldcatcodeofat

\newpage

\head {\S0 Introduction} \endhead  \resetall
\bigskip

We deal with the problem of the existence of a universal member in
${\frak K}_\lambda$.  For ${\frak K}$ a class of abelian groups, 
${\frak K}_\lambda$ is the class of $G \in {\frak K}$ of cardinality
$\lambda$; universal means every other member can be embedded into it.
We are concerned mainly with the class of reduced torsion free groups.  
Generally, on the history of the existence of universal members see 
Kojman-Shelah \cite{KjSh:409}.  From previous works, a natural
division of the possible cardinals for such problems is:
\bn
\ub{Case 0}:  $\lambda = \aleph_0$.  
\ub{Case 1}:  $\lambda = \lambda^{\aleph_0}$.
\ub{Case 2}:  $\aleph_0 < \lambda < 2^{\aleph_0}$
\ub{Case 3}:  $2^{\aleph_0} + \mu^+ < \lambda = \text{ cf}(\lambda) <
\mu^{\aleph_0}$.
\ub{Case 4}:  $2^{\aleph_0} + \mu^+ + \text{ cf}(\lambda) < \lambda <
\mu^{\aleph_0}$.
\ub{Case 5}:  $\lambda = \mu^+$, cf$(\mu) = \aleph_0,(\forall \chi < \mu)
(\chi^{\aleph_0} < \mu)$.
\ub{Case 6}:  cf$(\lambda) = \aleph_0,(\forall \chi < \lambda)(\chi^{\aleph_0}
< \lambda)$.
\ub{Subcase 6a}:  $\lambda$ is strong limit, cf$(\lambda) = \aleph_0$.
\ub{Subcase 6b}:  Case 6 but not 6a.
\bn
Our main interest will be in Case 3 for ${\frak K} = {\frak K}^{\text{rtf}}$, 
the class
of torsion free reduced abelian groups.  Note that divisible torsion free 
abelian groups of cardinality $\lambda$ are universal.  A second class 
is ${\frak K}^{\text{rs(p)}}$, the
class of reduced separable $p$-groups (see Definition \scite{2.3}(4), more in
Fuchs \cite{Fu}).  Kojman-Shelah \cite{KjSh:455} show that for ${\frak K} =
{\frak K}^{\text{rtf}},{\frak K}^{\text{rs(p)}}$ in Case 3 there is no 
universal member \ub{if} we restrict the possible embeddings to pure 
embeddings.  This stresses
that universality depends not only on the class of structures but also on the
kind of embeddings.  In \cite{Sh:456} we allow any embeddings, \ub{but}
restrict the class of abelian groups to $(< \lambda)$-stable ones.  In
\cite[\S1,\S5]{Sh:552} we allow any embedding and all $G \in {\frak K}_\lambda$
\ub{but} there is a further restriction on $\lambda$ related to the pcf
theory.  This restriction is weak in the following sense: it is not clear 
if there is any cardinal (in any possible universe of
set theory) not satisfying it.  We here prove the full theorem for $\lambda
> \beth_\omega$ with no further restrictions:
\mr
\item "{$(*)$}"   for $\lambda > \beth_\omega$ in Case 3, 
${\frak K} = {\frak K}^{\text{rtf}},{\frak K}^{\text{rs(p)}}$ there is no
universal member in ${\frak K}_\lambda$
\ermn
(where we define inductively $\beth_0 = \aleph_0,\beth_{n+1} = 2^{\beth_n},
\beth_\omega = \sum 2^{\beth_n}$ and generally $\beth_\alpha = \aleph_0 +
\dsize \sum_{\beta < \alpha} 2^{\beth_\beta}$).

\S1 deals with ${\frak K}^{\text{rtf}}$ using mainly type theory. 
In \S2, we apply combinatorial ideals whose definition has some built-in 
algebra and purely combinatorial ones to get results on 
${\frak K}^{\text{rs(p)}}$; there is more interaction between algebra and
combinatorics than in \cite{Sh:552}.  Similarly in \S3 we
work on the class of $\aleph_1$-free abelian groups.
\mn
What about the other cases?   Case 4 (like 3 but $\lambda$ singular)
for ${\frak K}^{\text{rtf}}_\lambda$ and pure embedding, was solved showing 
non-existence of universals in \cite{KjSh:455} if some weak pcf assumption 
holds and in \cite{Sh:552} this was done for embeddings under slightly 
stronger pcf assumptions.  For both assumptions, it is not clear if they 
may fail.  Note that the results on consistency of existence of
universals in this case cannot be attacked as long as more basic pcf 
problems remain open.
\mn
Concerning Case 5 - If we try to prove the consistency of the existence of
universals, it is natural first to prove the existence of the relevant club 
guessing; here we expect consistency results.  (Of course, consistently there 
is club guessing \nl
(by $\bar C = \langle C_\delta:\delta \in S \rangle,
S \subseteq \lambda$, otp$(C_\delta) = \mu$) and then there is no universal.)
Also we were first of all interested in the universal reduced torsion 
free groups under embeddings, but later we also looked into some of the 
other cases here.  See more in \cite{DjSh}.
\mn
Case 1 $(\lambda = \lambda^{\aleph_0})$.  By subsequent work there is a 
universal member of ${\frak K}^{\text{rtf}}_\lambda$, and 
(see Fuchs \cite{Fu}) in ${\frak K}^{\text{rs(p)}}_\lambda$ there is a 
universal member, but in ${\frak K}^{\aleph_1\text{-free}}_\lambda$ 
there is no universal member (see forthcoming work).
\mn
Case 0 $(\lambda = \aleph_0)$.  In ${\frak K}^{\text{rtf}}_\lambda$ there 
is no universal member (see above) and in ${\frak K}^{\text{rs(p)}}_\lambda$ 
there is a universal member (see Fuchs \cite{Fu}).
\mn
Case 2 $(\aleph_0 < \lambda < 2^{\aleph_0})$.  For ${\frak K}^{\text{rtf}}
_\lambda$ we prove here that there is no universal member (by \scite{1.2}), 
whereas for ${\frak K}^{\text{rs(p)}}_\lambda$ this is consistent with and 
independent of ZFC (see \cite[\S4]{Sh:550}). \nl
We also deal with Case 6 $((\forall \chi < \lambda)\chi^{\aleph_0} < \lambda,
\lambda > \text{ cf}(\lambda) = \aleph_0)$.  There is a universal member for 
${\frak K}^{\text{trf}}_\lambda$ and also for ${\frak K}^{\text{rs(p)}}
_\lambda$.  See \cite{DjSh}.
\mn
We thank two referees, Mirna Dzamonja and Norm Greenberg for many corrections.
\bn
\ub{Notation}:  The cardinality of a set $A$ is $|A|$, the cardinality of a
structure $G$ is $\|G\|$. \nl
${\Cal H}(\lambda^+)$ is the set of sets whose transitive closure has
cardinality $\le \lambda$ and $<^*_{\lambda^+}$ denotes a canonical well order
of ${\Cal H}(\lambda^+)$.

For an ideal $I$, we use $I^+$ to denote the family of subsets of Dom$(I)$
which are not in $I$.
\newpage

\head {\S1 Non-Existence of Universals Among Reduced Torsion Free Abelian
Groups} \endhead  \resetall
\bigskip

The first result (\scite{1.2}) deals with $\lambda$ satisfying
$\aleph_0 < \lambda < 2^{\aleph_0}$ and show the non-existence of universal
members in ${\frak K}^{\text{trf}}_\lambda$ which improves \cite{Sh:552}.  
The proof is straightforward by analyzing subgroups and comparing Bauer's 
types.  

Then we deal with $2^{\aleph_0} + \mu^+ < \lambda = \text{ cf}
(\lambda) < \mu^{\aleph_0}$.  We add witnesses to bar the way 
against ``undesirable" extensions (see \cite{DjSh} on classes of modules) 
which is a critical new point compared to \cite{Sh:552}.
\bigskip

\definition{\stag{1.1} Definition}  Let ${\frak K}^{\text{rtf}}$ 
denote the class 
of torsion free reduced abelian groups $G$ where torsion free means
$nx = 0,n \in \Bbb Z,x \in G \Rightarrow n = 0 \vee x = 0$ and reduced means 
$(\Bbb Q,+)$ cannot be embedded into $G$.  The subclass of $G \in
{\frak K}^{\text{rtf}}$ of cardinality $\lambda$ is denoted by 
${\frak K}^{\text{rtf}}_\lambda$.  Moreover, ${\frak K}^{\text{tf}}$ 
is the class of torsion free abelian groups.
\enddefinition
\bigskip

\proclaim{\stag{1.2} Claim}  1) If $\aleph_0 < \lambda < 2^{\aleph_0}$ 
\ub{then} ${\frak K}^{\text{rtf}}_\lambda$ has no universal member. \newline
2) Moreover, there is no member of ${\frak K}^{\text{rtf}}_\lambda$ 
universal for ${\frak K}^{\text{rtf}}_{\aleph_1}$.
\endproclaim
\bigskip

\demo{Proof}  Let $\bold P^*$ be the set of all primes and let 
$\{\bold Q_i:i < 2^{\aleph_0}\}$ be a family of infinite subsets of 
$\bold P^*$, pairwise with finite intersection.  Let $\rho_\alpha \in 
{}^\omega 2$ for $\alpha < \omega_1$ be pairwise distinct.  Let $H^*$ 
be the divisible torsion free abelian group with 
$\{x_\alpha:\alpha < \omega_1\}$ a maximal
independent subset.  For $i < 2^{\aleph_0}$ let $H^*_i$ be the subgroup
of $H^*$ generated by

$$
\align
\{x_\alpha:\alpha < \omega_1\} &\cup \{p^{-n}x_\alpha:p \in \bold P^* 
\backslash \bold Q_\alpha,\alpha < \omega_1 \text{ and } n < \omega\}\\
  &\cup \{p^{-n}(x_\alpha-x_\beta):\alpha,\beta < \omega_1 \text{ and }
p \in {\bold P}^* \text{ and} \\
  &\qquad \qquad \qquad \qquad \,\,\rho_\alpha \restriction p = 
\rho_\beta \restriction p \text{ and } n < \omega\}.
\endalign
$$
\medskip

\noindent
Clearly $H^*_i \in {\frak K}^{\text{rtf}}$ and 
$\|H^*_i\| = \aleph_1 \le \lambda$.
Let $G \in {\frak K}^{\text{rtf}}_\lambda$, and we shall prove 
that at most $\lambda$ of the groups $H^*_i$ are embeddable into $G$.

So assume $Y \subseteq 2^{\aleph_0},|Y| > \lambda$ and for $i \in Y$ we have
$h_i$ an embedding of $H^*_i$ into $G$ and we shall derive that $G$ is not
reduced; a contradiction.  We choose by induction on $n$ a set 
$\Gamma_n \subseteq {}^n \lambda$ and
pure abelian subgroups $G_\eta$ of $G$ for $\eta \in \Gamma_n$, as follows.
For $n=0$ we let $\Gamma_0 = \{<>\}$ and let $G_{<>} = G$.  For $n+1$, for
$\eta \in \Gamma_n$ such that $\|G_\eta\| > \aleph_0$ we let
$\Gamma_{n,\eta} = \{\eta \char 94 \langle \zeta \rangle:\zeta <
\|G_\eta\|\}$, and let $\bar G_\eta = \langle G_{\eta \char 94 \langle 
\zeta \rangle}:\zeta < \|G_\eta\| \rangle$ be an increasing continuous 
sequence of subgroups of $G_\eta$ of cardinality $< \|G_\eta\|$ with union
$G_\eta$ such that:
\medskip
\roster
\item "{$(*)$}"  for $\zeta < \|G_\eta\|$ we have \newline
$G_{\eta \char 94 \langle \zeta \rangle} = G_\eta \cap (\text{the Skolem Hull
of } G_{\eta \char 94 \langle \zeta \rangle} \text{ in } ({\Cal H}(\lambda^+),
\in,<^*_{\lambda^+},G_\eta))$.
\endroster
\medskip

\noindent
Let $\Gamma_{n+1} = \{\eta \char 94 \langle \zeta \rangle:\eta \in
\Gamma_n,\|G_\eta\| > \aleph_0 \text{ with } \zeta < \|G_\eta\|\}$ and
$\Gamma = \dsize \bigcup_{n < \omega} \Gamma_n$.

For each $i \in Y$, let $\eta = \eta_i \in \Gamma$ be such that:
\medskip
\roster
\item "{$(a)$}"  $\{\alpha < \omega_1:h_i(x_\alpha) \in G_{\eta_i}\}$ is
uncountable
\sn
\item "{$(b)$}"  under $(a)$, the cardinality of $G_{\eta_i}$ is minimal.
\endroster
\medskip

\noindent
Clearly $\eta_i$ is well defined as $(a)$ holds for $\eta = \langle \rangle$
and clearly $G_{\eta_i}$ is uncountable.
Let $X_i = \{\alpha < \omega_1:h_i(x_\alpha) \in G_{\eta_i}\}$, and let
$\beta_i < \omega_1$ be minimal such that \nl
$\{\rho_\alpha:\alpha \in \beta_i \cap X_i\}$ is a dense subset of 
$\{\rho_\alpha:\alpha \in X_i\}$.  Let $\zeta_i < \|G_{\eta_i}\|$ 
be the minimal $\zeta$ such that 
$\{h_i(x_\alpha):\alpha \in \beta_i \cap X_i\} \subseteq 
G_{\eta \char 94 \langle \zeta \rangle}$.  Now the set \nl
$X'_i = \{\alpha < \omega_1:h_i(x_\alpha) \in 
G_{\eta_i \char 94 \langle \zeta_i \rangle}\}$ is countable, and hence we can 
find $\alpha_i \in X_i \backslash X'_i$.
\medskip

Now the number of possible sequences $\langle \eta_i,\beta_i,\zeta_i,
\alpha_i,h_i(x_{\alpha_i}) \rangle$ is at most $|{}^{\omega >}\lambda| 
\times \aleph_1 \times \lambda \times \aleph_1 \times \lambda$ 
(as $\Gamma \subseteq {}^{\omega >} \lambda$). 
So for some $\langle \eta,\beta,\zeta,\alpha,y \rangle$ and $i_0 < i_1$ from
$Y$ we have (for $\ell=0,1$)

$$
\eta_{i_\ell} = \eta,\beta_{i_\ell} = \beta,\zeta_{i_\ell} = \zeta,
\alpha_{i_\ell} = \alpha,h_{i_\ell}(x_{\alpha_\ell}) = y.
$$
\medskip

\noindent
Now as $h_{i_\ell}$ embeds $H^*_{i_\ell}$ into $G$ and $h_{i_\ell}(x_\alpha)
=y$ necessarily
\medskip
\roster
\item "{$(*)$}"  if $p \in \bold P^* \backslash \bold Q_{i_\ell}$ and 
$n < \omega$ then in $G,p^{-n}$ divides $y$.
\endroster
\medskip

\noindent
So this holds for every $p \in (\bold P^* \backslash \bold Q_{i_0}) 
\cup (\bold P^* \backslash \bold Q_{i_1}) = \bold P^* \backslash 
(\bold Q_{i_0} \cap \bold Q_{i_1})$.

Now $\bold Q_{i_0} \cap \bold Q_{i_1}$ is finite so let $p^* \in \bold P^*$ 
be above its supremum.  As $\{\gamma:\gamma \in X'_{i_0}\}$ is a 
dense subset of $\{\rho_\alpha:\alpha \in X_{i_0}\}$, there is 
$\gamma \in X'_{i_0}$ such that $\rho_\gamma \restriction
p^* = \rho_\alpha \restriction p^* (= \rho_{\alpha_{i_0}} \restriction p^*)$.
Let $h_{i_0}(x_\gamma) = y^*$, it is in $G_{\eta \char 94 \langle \zeta 
\rangle}$.
\smallskip
So in $({\Cal H}(\lambda^+),\in,<^*_{\lambda^+},G_\eta)$, the following
formula is satisfied

$$
\align
\varphi(y,y^*) = ``\text{in } G_\eta,\,&y \text{ is divisible by } p^n
\text{ when } p \in \bold P^* \and p \ge p^* \and n < \omega \\
  &\text{and } y-y^* \text{ is divisible by } p^n \text{ when}\\
  &p \in \bold P^* \and p < p^* \and n < \omega".
\endalign
$$
\medskip

\noindent
Hence by $(*)$, i.e. by the choice of 
$\langle G_{\eta \char 94 \langle \xi \rangle}:\xi <
\|G_\eta\| \rangle$ necessarily for some $y' \in 
G_{\eta \char 94 \langle \zeta \rangle}$ we have $\varphi(y',y^*)$.  Now
$y \ne y'$ as $y' \in G_{\eta \char 94 \langle \zeta \rangle},y \notin
G_{\eta \char 94 \langle \zeta \rangle}$.  Also $y-y'$ is divisible by $p^n$
for $p \in \bold P^*,n < \omega$. \newline
[Why?  If $p \ge p^*$ because both $y$ and $y'$ are divisible by $p^n$ 
and if $p < p^*$ because $y-y' = (y-y^*) - (y'-y^*)$ and both $y-y^*$ 
and $y'-y^*$ are divisible by $p^n$.] \newline
As $G$ is torsion free, the pure closure in $G$ of $\langle \{y-y'\} \rangle
_G$ is isomorphic to $\Bbb Q$, a contradiction to ``$G$ is reduced".
\hfill$\square_{\scite{1.2}}$
\enddemo
\bigskip

\definition{\stag{1.3} Definition}  1) Let $\bold P^*$ be the set of primes.
\newline
2) For $G \in {\frak K}^{\text{rtf}}$ and $x \in G \backslash \{0\}$ let
\medskip
\roster
\item "{$(a)$}"  $\bold P(x,G) = 
\{p \in \bold P^*:x \in \dsize \bigcap_{n < \omega} p^n G$, \newline

$\qquad \qquad \qquad \qquad$ equivalently $x$ is divisible by $p^n$ \newline

$\qquad \qquad \qquad \qquad$ in $G$ for every $n < \omega\}$
\smallskip
\noindent
\item "{$(b)$}"  $\bold P^-(x,G) = 
\{p:p \in \bold P^*$, but $p \notin \bold P(x,G)$ \newline

$\qquad \qquad \qquad$ and there is $y \in G \backslash \{0\}$ such that
\newline

$\qquad \qquad \qquad \bold P^* \backslash \{p\} \subseteq \bold P(y,G)$
and $p \in \bold P(x-y,G)\}$.
\endroster
\medskip

\noindent
3) $G \in {\frak K}^{\text{rtf}}$ is called full if: for every 
$x \in G \backslash \{0\}$ we have $\bold P^* = \bold P(x,G) 
\cup \bold P^-(x,G)$. \newline
4) The class of full $G \in {\frak K}^{\text{rtf}}$ is 
called ${\frak K}^{\text{stf}}$ and ${\frak K}^{\text{stf}}_\lambda = 
{\frak K}^{\text{stf}} \cap {\frak K}^{\text{rtf}}_\lambda$,
(why $s$?  as the successor of $r$ in the alphabet).
\enddefinition
\bigskip

\demo{\stag{1.4} Fact}  1) If $G \in {\frak K}^{\text{rtf}}$, then for 
any $x \in G$ the sets $\bold P(x,G)$ and $\bold P^-(x,G)$ are 
disjoint subsets of $\bold P^*$. \newline
2) If $G_2$ is an extension of $G_1$, both in ${\frak K}^{\text{rtf}}$ and
$x \in G_1 \backslash \{0\}$ \ub{then}
\medskip
\roster
\item "{$(a)$}"  $\bold P(x,G_1) \subseteq \bold P(x,G_2)$, with equality
if $G_1$ is a pure subgroup of $G_2$
\sn
\item "{$(b)$}"  $\bold P^-(x,G_1) \subseteq \bold P^-(x,G_2)$.
\endroster
\medskip

\noindent
3)  For every $G \in {\frak K}^{\text{rtf}}$ there is a $G'$ such that
\medskip
\roster
\item "{$(a)$}"  $G'$ is full, $G' \in {\frak K}^{\text{rtf}}$
\sn
\item "{$(b)$}"  $G$ is a pure subgroup of $G',\|G'\| = \|G\|$.
\endroster
\enddemo
\bigskip

\demo{Proof}  1),2)  Trivial. \newline
3) It suffices to show
\medskip
\roster
\item "{$(*)$}"  if $G \in {\frak K}^{\text{rtf}}$ 
and $x \in G \backslash \{0\}$,
and $p \in \bold P^* \backslash \bold P(x,G)$ \underbar{then} for some pure
extension $G'$ of $G$ with rk$(G/G')=1$ we have: 
$p \in \bold P^-(x,G')$.
\endroster
\medskip

\noindent
For proving $(*)$ for a given $G,x$ let $\hat G$ be the divisible hull of
$G$ and let \nl
$G_0 = \{y \in \hat G:\text{ for some } n>0,p^n y \in G\}$, \nl
$G_1 = \{y \in \hat G:\text{ for some } b \in \Bbb Z,b > 0 
\text{ not divisible by } p \text{ we have by } \in G\}$.  
Clearly $G=G_0 \cap G_1$.  
We define the following subsets of $\hat G \times \Bbb Q$:

$$
H_0 = \{(y,0):y \in G\} \text{ (so } G \text{ is isomorphic to } H_0)
$$

$$
H_1 = \{(p^nbx,p^nb):b,n \in \Bbb Z\}
$$

$$
H_2 = \{(0,c_1/c_2):c_1,c_2 \in \Bbb Z \text{ and } c_2 \text{ not divisible
by } p\}.
$$
\mn
Easily all three are additive subgroups of $\hat G \times \Bbb Q$ and 
$H_2 \cong \Bbb Z_{(p)}$.  Let $G' = H_0 + H_1 + H_2$, a subgroup of
$\hat G \times \Bbb Q$. \newline
We claim that $G' \cap (\hat G \times \{0\}) = H_0$.  The inclusion
$\supseteq$ should be clear.  For the other direction let $z \in G' \cap
(\hat G \times \{0\})$; as $z \in G'$ there are $(y,0) \in
H_0$, (so $y \in G$), $(p^nbx,p^nb) \in H_1$ (so $b \in \Bbb Z,n \in \Bbb Z$ 
and $x \in G$ is the constant from $(*)$) and 
$(0,c_1/c_2) \in H_2$ (so $c_1,c_2 \in \Bbb Z$ and $p$ does 
not divide $c_2$) and integers $a_0,a_1,a_2$ such that
$z = a_0(y,0) + a_1(p^nbx,p^nb) + a_2(0,c_1/c_2)$ which means
$z = (a_0y + a_1p^nbx,a_1p^nb + a_2c_1/c_2)$.
\mn
As $z \in \hat G \times \{0\}$ clearly $a_1p^nb+a_2c_1/c_2 =0$, so as $p$
does not divide $c_2$, necessarily $a_1p^nb$ is an integer, hence
$a_1p^nbx \in G$, hence $a_0y+a_1p^nbx \in G$ and hence 
$z \in G \times \{0\} = H_0$ as required.  \nl
It is easy to check now that $H_0$ is a pure subgroup of $G'$. \nl
Also letting $y^* = (0,-1)$ clearly $(x,0) - y^*$ is divisible by
$p^k$ for every $k < \omega$ (as $(p^k x,p^k) \in H_1 \subseteq G'$ for every
$k \in \Bbb Z$) and $y^*$ is
divisible by any integer $b$ when $b$ is not divisible by $p$ (as
$\frac 1b\,\,y^* = (0,-1/b) \in H_2 \subseteq G'$).
\mn
Identifying $y \in G$ with $(y,0) \in G$ we are done: $G'$ is as
required in $(*)$, with $y^*$ witnessing ``$p \in \bold P^-(x,G')$". 
\hfill$\square_{\scite{1.4}}$
\enddemo
\bigskip

\proclaim{\stag{1.5} Claim}  If $G_1 \in {\frak K}^{\text{\rm rtf\/}}$ 
is full and $G_2 \in {\frak K}^{\text{\rm rtf\/}}$ and $h$ is an 
embedding of $G_1$ into $G_2$ then:

$$
\text{for } x \in G_1 \backslash \{0\},\bold P(x,G_1) = \bold P(h(x),G_2).
$$
\endproclaim
\bigskip

\demo{Proof}  Without loss of generality $h$ is the identity, now reflect
using \scite{1.4}(1), \scite{1.4}(2) and the definition of full.
\hfill$\square_{\scite{1.5}}$
\enddemo
\bigskip

\demo{\stag{1.6} Conclusion}  Assume
\medskip
\roster
\item "{$(*)$}"  $2^{\aleph_0} < \mu^+ < \lambda = \text{ cf}(\lambda) <
\mu^{\aleph_0}$.
\endroster
\medskip

\noindent
Then there is no universal member in ${\frak K}^{\text{rtf}}_\lambda$.
\enddemo
\bigskip

\demo{Proof}  Let $S \subseteq \{\delta < \lambda:\text{cf}(\delta) =
\aleph_0 \text{ and } \omega^2 \text{ divides } \delta\}$ be stationary 
and $\bar \eta = \langle \eta_\delta:\delta \in S \rangle$ where each
$\eta_\delta$ is an increasing $\omega$-sequence
of ordinals $< \delta$ with limit $\delta$ such that $\eta_\delta(n)-n$ is
divisible by $\omega$; so $\delta_1 \ne \delta_2 \Rightarrow
\text{ Rang}(\eta_{\delta_1}) \cap \text{ Rang}(\eta_{\delta_2})$ is finite.
Let $\{p^*_n:n < \omega\}$ list the primes in the increasing order.  Let
$G^0_{\bar \eta}$ be the abelian group generated by $\{x_\alpha:\alpha <
\lambda\} \cup \{y_\delta:\delta \in S\} \cup \{z_{\delta,n,\ell}:n,\ell
< \omega\} \cup \{z_{\delta,n,m,\ell}:n,\ell < \omega,m \in \omega
\backslash \{n\}\}$ freely except for the equations

$$
p^*_nz_{\delta,n,\ell +1} = z_{\delta,n,\ell} \qquad
y_\delta - x_{\eta_\delta(n)} = z_{\delta,n,0}.
$$

$$
p^*_m z_{\delta,n,m,\ell +1} = z_{\delta,n,m,\ell +1},x_{\eta_\delta} =
z_{\delta,n,m,0}.
$$
\medskip

\noindent
We can check that $G^0_{\bar \eta} \in {\frak K}^{\text{rtf}}_\lambda$.
\medskip

Let $G_{\bar \eta} \in {\frak K}^{\text{rtf}}_\lambda$ be a pure extension of
$G^0_{\bar \eta}$ which is full (one exists by \scite{1.4}(3)).  So
\medskip
\roster
\item "{$(*)$}"  if $h$ embeds $G_{\bar \eta}$ into $G \in 
{\frak K}^{\text{rtf}}_\lambda$ then \newline
$x \in G_{\bar \eta} \backslash \{0\} \Rightarrow \bold P(x,G_{\bar \eta})
= \bold P(h(x),G)$.
\endroster
\medskip

\noindent
Hence the proof in \cite{KjSh:455} works. \hfill$\square_{\scite{1.6}}$
\enddemo
\bigskip

\remark{\stag{1.7} Remark}  1) Similarly the results on $\lambda$ singular
(i.e. Case 4) in \cite{KjSh:455}, hold for embedding (rather than pure
embedding). \nl
2) What about Case 5?  If there is a family ${\Cal P} \subseteq \{C \subseteq
\mu^+:\text{otp}(C) = \mu\}$ which guesses clubs (i.e. every club
$E$ of $\mu^+$ contains one of them), the result holds. \nl
3) On $\aleph_0 \le \lambda < 2^{\aleph_0}$ see also in \scite{3.15}.
\endremark
\newpage

\head {\S2 The existence of universals for separable reduced abelian
$p$-groups} \endhead  \resetall
\bigskip

We here eliminate the very weak $pcf$ assumption from the theorem of
``no universal in ${\frak K}^{\text{rs(p)}}_\lambda$" when 
$\lambda > \beth_\omega$.  Note that ${\frak K}^{\text{rs(p)}}$ is 
defined in \scite{2.3}(4).

In the first section we have eliminated the very weak $pcf$ assumptions for
the theorem concerning ${\frak K}^{\text{rtf}}_\lambda$ (though the 
$\lambda = \text{ cf}(\lambda) > \mu^+$ remains, i.e. we assume we are in
Case 3).  
This was done using the ``infinitely many primes", so in the language of
e.g. \cite{KjSh:455} the invariant refers to one element $x$.  This cannot
be generalized to ${\frak K}^{\text{rs(p)}}_\lambda$.  
However, in \cite[\S5]{Sh:552}
we use an invariant on e.g. suitable groups and related stronger
``combinatorial" ideals.  We continue this, using combinatorial ideals
closer to the algebraic ones to show that the algebraic is non-trivial.

We rely on the ``GCH" proved in \cite{Sh:460} hence the condition ``$\lambda >
\beth_\omega$" is used.
\bigskip

\definition{\stag{2.1} Definition}  1) For $\bar \lambda = 
\langle \lambda_\ell:\ell < \omega \rangle$ and $\bar t = \langle t_\ell:
\ell < \omega \rangle$ (with $1 < t_\ell < \omega$) we define 
$J^4_{\bar t,\bar \lambda}$.

It is a family of subsets $A$ of $\dsize \prod_{\ell < \omega}
[\lambda_\ell]^{t_\ell}$ such that: \nl
\medskip
\roster
\item "{$(*)$}"  for every large enough $\ell < \omega$, for every 
$B \in [\lambda_\ell]^{\aleph_0}$ for some $k \in (\ell,\omega)$ we cannot
find 
$$
\langle \nu_\eta:\eta \in \dsize \prod_{i \in [\ell,k)} [\omega]^{t_i}
\rangle
$$
\noindent
such that
{\roster
\itemitem{ (a) }  $\nu_\eta \in A$
\smallskip
\noindent
\itemitem{ (b) }  if $\eta_1,\eta_2 \in \dsize \prod_{i \in [\ell,k)} 
[\omega]^{t_i},\ell \le m \le k$ and $\eta_1 \restriction [\ell,m) = \eta_2 
\restriction [\ell,m)$ \ub{then} $\nu_{\eta_1} \restriction m = \nu_{\eta_2}
\restriction m$; hence \nl
$\nu_{\eta_1} \restriction \ell = \nu_{\eta_2} \restriction
\ell$ for $\eta_1,\eta_2 \in \dsize \prod_{i \in [\ell,k)}[\omega]^{t_i}$
\smallskip
\noindent
\itemitem{ (c) }  if $\eta_0 \in \dsize \prod_{i \in [\ell,k)}[\omega]^{t_i}$
and $\ell \le m < k$ \ub{then} for some 
$E \in [\lambda_m]^{\aleph_0}$ we have
$$
\qquad [E]^{t_m} = \{\nu_\eta(m):\eta \in 
\dsize \prod_{i \in [\ell,k)}[\omega]^{t_i} \text{ and } 
\eta \restriction m = \eta_0 \restriction m\}
$$
\noindent
and $m = \ell \Rightarrow E=B$.
\endroster}
\endroster
\medskip

\noindent
2) Let 

$$
\align
J^4_{\bar t,\bar \lambda,< \theta} = \bigl\{ A \subseteq \dsize \prod
_{\ell < \omega}[\lambda_\ell]^{t_\ell}:&\text{ for some } \alpha < \theta
\text{ and } A_\beta \in J^4_{\bar t,\bar \lambda} \text{ for } 
\beta < \alpha \\
  &\text{ we have }
A \subseteq \dsize \bigcup_{\beta < \alpha} A_\beta \bigr\}.
\endalign
$$
\medskip

\noindent
When $\theta = \kappa^+$, we may write $\kappa$ instead of $< \theta$.
\enddefinition
\bigskip

\demo{\stag{2.2} Fact}  1) $J^4_{\bar t,\bar \lambda,\theta}$ is a
$\theta$-complete ideal. \newline
2) If $\ell < \omega \Rightarrow \lambda_\ell > \beth_{t_\ell-1}(\theta)$ 
then the ideal $J^4_{\bar t,\bar \lambda,\theta}$ is proper (where
$\beth_0(\theta) = \theta,\beth_{n+1}(\theta) = 2^{\beth_n(\theta)}$,
and for $\alpha$ limit $\beth_\alpha(\theta) = \theta + 
\dsize \sum_{\beta < \alpha} 2^{\beth_\beta(\theta)}$.
\enddemo
\bigskip

\demo{Proof}  1) Trivial. \newline
2) Let for $\ell < \omega$

$$
\align
ERI^{n_\ell}_{\lambda_\ell} = \bigl\{A \subseteq 
[\lambda_\ell]^{t_\ell}:&\text{ for some } F:
[\lambda_\ell]^{t_\ell} \rightarrow \theta \text{ there is no } 
B \in [\lambda_\ell]^{\aleph_0} \\
  &\text{ such that } 
F \restriction [B]^{t_\ell} \text{ is constant } \text{ and } 
[B]^{t_\ell} \subseteq A \bigr\}.
\endalign
$$
\medskip

\noindent
So this is a $\theta^+$-complete ideal.  It is non-trivial by 
Erd\"os-Rado theorem (we use it similarly in \cite[\S1]{Sh:620}).  Now we
shall prove that the ideal $J^4_{\bar t,\bar \lambda,\theta}$ is proper. 
So assume
$\dsize \prod_{\ell < \omega}[\lambda_\ell]^{t_\ell} = \dsize \bigcup
_{i < \theta} X_i$ and $X_i \in J^4_{\bar t,\bar \lambda}$ for each
$i < \theta$ and we shall get a contradiction.  Let

$$
X^+_i = \bigl\{ \eta \in \dsize \prod_{\ell < \omega}
[\lambda_\ell]^{t_\ell}:\text{ for every } \ell < \omega \text{ for some }
\eta' \in X_i \text{ we have } \eta \restriction \ell = \eta' 
\restriction \ell \bigr\}.
$$
\medskip

\noindent
(i.e. the closure of $X_i$).  So $X^+_i \subseteq \dsize \prod_{\ell < \omega}
[\lambda_\ell]^{t_\ell} = \dsize \prod_{\ell < \omega} \text{ Dom}
(ERI^{t_\ell}_{\lambda_\ell})$ is closed, and those ideals are 
$\theta^+$-complete
and $\dsize \prod_{\ell < \omega} \text{ Dom}(ERI^{t_\ell}_{\lambda_i}) =
\dsize \bigcup_{i < \theta} X^+_i$.  Hence (Rubin-Shelah \cite{RuSh:117},
\cite[Ch.XI,3.5(2)]{Sh:f} with $H_\alpha = X^+_i$) we can find $T$ 
such that:
\medskip
\roster
\item "{$(a)$}"  $T \subseteq \dsize \bigcup_{m < \omega}
\dsize \prod_{\ell < m}[\lambda_\ell]^{t_\ell}$
\sn
\item "{$(b)$}"  $T$ is closed under initial segments
\sn
\item "{$(c)$}"  $<> \in T$
\sn
\item "{$(d)$}"  if $\nu \in T$ and $\ell g(\nu) = \ell$ then \newline
$\{u \in [\lambda_\ell]^{t_\ell}:\nu \char 94 \langle u \rangle \in
T\} \in (ERI^{t_\ell}_{\lambda_\ell})^+$
\sn
\item "{$(e)$}"  for some $i < \theta$, lim$(T) \subseteq X^+_i$.
\ermn
(Here, lim$(T) = \{\nu \in \dsize \prod_{\ell < \omega}[\lambda_\ell]
^{t_\ell}:(\forall m < \omega)\nu \restriction m \in T\})$. \nl
By the definition of the ideal we get more than required (for every $k$
in place of ``some $k$" in $(*)$ of Definition \scite{2.1}).
\hfill$\square_{\scite{2.2}}$
\enddemo
\bigskip

\remark{Remark}  So we could have used the stronger ideal defined implicitly
in \scite{2.2}, i.e. $J^5_{\bar t,\bar \lambda,\theta} = \{X \subseteq
\dsize \prod_{\ell < \omega} \lambda_\ell:\text{ we can find } \alpha <
\theta \text{ and } X_i \subseteq X \text{ for } i < \alpha \text{ such that }
X = \dbcu_{i < \alpha} X_i \text{ and for each } i \text{ and } T$ satisfying
clauses $(a)-(d)$ from the proof of \scite{2.2} there is $T' \subseteq T$
satisfying clauses $(a)-(d)$ such that lim$(T)$ is disjoint to the closure of
$X_i\}$.  Of course, we can also replace ERI$^{t_\ell}_{\lambda_\ell}$ by 
various variants.
\endremark
\bigskip

\noindent
We recall from \cite[5.1]{Sh:552}
\definition{\stag{2.3} Definition}  (\cite[5.1]{Sh:552})  1) For $\bar \mu = 
\langle \mu_n:n < \omega \rangle$ let $B_{\bar \mu}$ be the following direct
sum of cyclic $p$-groups.  Let $K^n_\alpha$ be a cyclic group of order
$p^{n+1}$ generated by $x^n_\alpha$ and let $B_{\mu_n} = 
\dsize \oplus_{\alpha < \mu_n} K^n_\alpha$ and $B_{\bar \mu} = 
\dsize \oplus_{n < \omega} B_{\mu_n}$, i.e. $B_{\bar \mu}$ is the 
abelian group generated by $\{x^n_\alpha:n < \omega,\alpha < \mu_n\}$ 
freely except that $p^{n+1}\,x^n_\alpha=0$.
\sn
Moreover, let $B_{\bar \mu \restriction n} = 
\oplus \{K^m_\alpha:\alpha < \mu_m,m < n\} \subseteq B_{\bar \mu}$ \newline
(these groups are in ${\frak K}^{\text{rs(p)}}
_{\le \underset n {}\to \sum \mu_n}$).
\newline
Let $\hat B_{\bar \mu}$ be the $p$-torsion completion of $B_{\bar \mu}$ 
(i.e. completion under the norm \nl
$\|x\| = \text{ min}\{2^{-n}:p^n \text{ divides } x\}$ but putting only the
torsion elements, see Fuchs \cite{Fu}.  Note that $\hat B_{\bar \mu}$ is the
torsion part of the $p$-adic completion of $B_{\bar \mu}$). \newline
2) Let $I^1_{\bar \mu,<\theta} = I^1_{\bar \mu,\theta}[p]$ be the ideal 
on $\hat B_{\bar \mu}$ (depending on the choice of $\langle K^n_\alpha:
\alpha < \mu_n;n < \omega \rangle$ or actually $\langle 
B_{\bar \mu \restriction n}:n < \omega \rangle$) consisting of unions of 
$< \theta$ members of $I^0_{\bar \mu}$, where

$$
I^0_{\bar \mu} = I^0_{\bar \mu,\theta}[p] = \bigl\{ A \subseteq 
\hat B_{\bar \mu}:\text{for every large enough } n, \text{ we have }
  c \ell_{\hat B_{\bar \mu}}(\langle A \rangle_{\hat B_{\bar \mu}}) \cap
B_{\bar \mu} \subseteq B_{\bar \mu \restriction n} \bigr\}
$$
\mn
($c \ell_{\hat B_{\bar \mu}}$ is defined in 3) below). \nl
When $\theta = \kappa^+$ instead of $< \theta$ we may write $\kappa$.  
If $\mu_n = \mu$, we may write $\mu$ instead of $\bar \mu$. \nl
3) For $X \subseteq \hat B_{\bar \mu}$, recall

$$
c \ell_{\hat B_{\bar \mu}}(X) = \bigl\{ x:(\forall n)(\exists y \in X)
(x-y \in p^n \hat B_{\bar \mu}) \bigr\}.
$$
\mn
4) Let ${\frak K}^{\text{rs(p)}}$ be the family of pure subgroups of some
$\hat B_{\bar \mu}$. \nl
5) If $p$ is not clear from the context, we may write
$B^p_{\bar \mu},\hat B^p_{\bar \mu}$, etc.
\enddefinition
\bigskip

\proclaim{\stag{2.4} Claim}  Assume $\bar \mu = \langle \mu_n:n < \omega
\rangle,\bar t = \langle t_\ell:\ell < \omega \rangle,t_\ell = p$ and
the ideal $J^4_{\bar t,\bar \mu,\theta}$ is proper 
(so $\mu_n \ge \beth_{p-1}(\theta)^+$ is enough by \scite{2.2}(2)).  
\ub{Then} the ideal $I^1_{\bar \mu,\theta}$ is proper 
(and $I^1_{\bar \mu,\theta}$ is a $\theta^+$-complete ideal).
\endproclaim
\bigskip

\demo{Proof}  We define a function $h$ from $\dsize \prod_{\ell < \omega}
[\lambda_\ell]^{t_\ell}$ into $\hat B_{\bar \mu}$.  We let

$$
h(\eta) = \Sigma\{ p^n\,x^n_\beta:\beta \in \eta(n) \text{ and }
n < \omega\} \in \hat B_{\bar \mu}[p].
$$
\mn
Clearly $h$ is one to one and it suffices to prove
\medskip
\roster
\item "{$(*)$}"  if $X \in (J^4_{\bar t,\bar \mu,\theta})^+$ then 
$h''(X)$ belongs to $(I^1_{\bar \mu,\theta})^+$.
\endroster
\medskip

\noindent
So assume $X \in (J^4_{\bar t,\bar \lambda,\theta})^+$ is given and suppose
toward contradiction that $h''(X) \in I^1_{\bar \mu,\theta}$.  So we can
find $\langle Y_i:i < \theta \rangle$ such that for such $i < \theta$ we have
$Y_i \in I^0_{\bar \mu}$ and $h(X) \subseteq \dbcu_{i < \theta} Y_i$.  Let
$X_i = h^{-1}(Y_i)$.  So $h(X_i) \subseteq Y_i \in I^0_{\bar \mu}$ and hence
$h(X_i) \in I^0_{\bar \mu}$, but as $J^4_{\bar t,\bar \lambda,\theta}$ is
$\theta^+$-complete and $X \in (J^4_{\bar t,\bar \lambda,\theta})^+$ 
necessarily for some $i < \theta,X_i \in 
(J^4_{\bar t,\bar \lambda,\theta})^+$, so without loss of generality 
$h''(X) \in I^0_{\bar \mu}$.  By the definition of 
$I^0_{\bar \mu}$, for some $n(*) < \omega$ we have
\mr
\item "{$(*)$}"  $B_{\bar \mu} \cap c \ell(h''(X)) \subseteq B_{\bar \mu
\restriction n(*)}$.
\ermn
On the other hand, as $X \in (J^4_{\bar t,\bar \mu,\theta})^+$, it is
$\notin J^4_{\bar t,\bar \mu}$ so from definition \scite{2.1}(1) of
$J^4_{\bar t,\bar \mu}$ we can find $\langle B_n:n \in \Gamma \rangle$ 
such that:
\medskip
\roster
\item "{$(a)$}"  $\Gamma \in [\omega]^{\aleph_0}$ and $B_n \in [\lambda_n]
^{\aleph_0}$
\sn
\item "{$(b)$}"  for $n \in \Gamma$, for every $k \in (n,\omega)$ we can find
\newline
$\langle \nu^{n,k}_\eta:\eta \in \dsize \prod_{\ell \in [n,k)}
[\omega]^{t_\ell} \rangle$ as in (a)-(c) of Definition \scite{2.1}(1) with
$n,B_n,k$ here standing for $\ell,B,k$ there.
\endroster
\medskip

\noindent
For $m \in (n,k]$ and $\eta \in \dsize \prod_{\ell \in [n,m)}
[\omega]^{t_\ell}$ we let $\nu^{n,k}_\eta$ be $\nu^{n,k}_{\eta_1} 
\restriction m$ whenever \newline
$\eta \triangleleft \eta_1 \in \dsize \prod_{\ell \in [n,k)}
[\omega]^{t_\ell}$ (by clause (b) in $(*)$ of \scite{2.1} it is well 
defined).   Fix $n \in \Gamma$ and $k \in [n,\omega)$ for awhile.
Let $A_\eta = A^{n,k}_\eta \in [\lambda_m]^{\aleph_0}$ be such that 
$\{\nu^{n,k}_{\eta \char 94 \langle u \rangle}(m):u \in [\omega]^{t_m}\} = 
[A_\eta]^{t_m}$ and without loss of generality (otp stands for ``the order
type")
\medskip
\roster
\item "{$(*)$}"  otp$(A_\eta) = \omega$ and $\nu^{n,k}_{\eta \char 94 
\langle u \rangle}(m) = \text{ OP}_{A_\eta,\omega}(u)$
\endroster
\medskip

\noindent
(where $OP_{A_\eta,\omega}(i) = \alpha$ iff $i = \text{ otp}
(A_\eta \cap \alpha)$).
\nl
Now for $m \in (n,k]$ and $\eta \in \dsize \prod_{\ell \in [n,m)}[w]^{t_\ell}$
we define

$$
y_\eta = y^{n,k}_\eta = 
\dsize \sum \bigl\{ h(\nu^{n,k}_\rho):\eta \trianglelefteq \rho \in
\dsize \prod_{\ell \in [n,k)}[\omega]^{t_\ell} \text{ and }
  (\forall \ell)[\ell g(\eta) \le \ell < k \rightarrow \rho(\ell)
\subseteq [0,t_\ell] \bigr\}
$$
\mn
where $\trianglelefteq$ denotes being an initial segment.  
So $y_\eta \in \hat B_{\bar \mu}$ and we shall prove by downward induction
on $m \in (n,k]$ that for every $\eta \in \dsize \prod_{\ell \in [n,m)}
[w]^{t_\ell}$ we have
\mr
\item "{$\boxtimes$}"  $y_\eta = \bigl( \dsize \prod^{k-1}_{\ell=m}
(t_\ell +1) \bigr) \times \bigl( \dsize \sum_{\ell < m} \,\,
\dsize \sum_{\alpha \in \nu^{n,k}_\eta(\ell)}p^\ell \, x^\ell_\alpha \bigr)
\text{ mod } p^k \, \hat B_{\bar \mu}$.
\ermn
\ub{Case 1}:  $m = k$.

In this case the product $\dsize \prod^{k-1}_{\ell=m}(t_\ell +1)$ is just 1,
so the equation says \nl
$y_\eta = \dsize \sum_{\ell < m} \,\, \dsize \sum_{\alpha \in \nu^{n,k}_\eta
(\ell)} p^\ell \, x^\ell_\alpha \text{ mod } p^k \, \hat B_{\bar \mu}$. \nl
Now the expression for $y_\eta$ is

$$
\align
\dsize \sum \{ h(\nu^{n,k}_\rho):&\eta \trianglelefteq \rho \in
\dsize \prod_{\ell \in [n,k)}[\omega]^{t_\ell} \text{ and } 
(\forall \ell)[m \le \ell < k \Rightarrow \rho(\ell) \subseteq 
[0,t_\ell]] \} \\
  &= h(\nu^{n,k}_\eta) = \dsize \sum_{\ell < \omega} \,\,
\dsize \sum_{\alpha \in \nu^{n,k}_\eta(\ell)} p^\ell \, x^\ell_\alpha \\
  &= \dsize \sum_{\ell < m} \,\, 
\dsize \sum_{\alpha \in \nu^{n,k}_\eta(\ell)} p^\ell \, x^\ell_\alpha
 + p^k(\dsize \sum_{\ell \in [k,\omega)}\,\, \dsize \sum_{\alpha
\in \nu^{n,k}_\eta(\ell)} p^{\ell -k} x^\ell_\alpha)
\endalign
$$
\mn
so the equality is trivial.
\bn
\ub{Case 2}:  $n < m < k$.

Here (with equalities in the equation being in $\hat B_{\bar \mu}$, modulo 
$p^k \, \hat B_{\bar \mu}$), we have: 

$$
\align
y_\eta &= \\
  &\qquad \qquad \qquad \qquad \qquad \qquad \text{ [by the definition of } 
y_\eta,y_{\eta \char 94 \langle u \rangle}] \\
   &= \dsize \sum \{y_{\eta \char 94 \langle u \rangle}:u \in [\{0,\dotsc,
t_m\}]^{t_m}\} = \\
  &\qquad \qquad \qquad \qquad \qquad \qquad \text{[by the induction
hypothesis]} \\
&= \dsize \sum \{( \dsize \prod_{\ell = m+1}^{k-1}(t_\ell + 1))
(\dsize \sum_{\ell < m +1} \,\,\, \dsize \sum_{\alpha \in \nu^{n,k}
_{\eta \char 94 \langle u \rangle}(\ell)} p^\ell \, x^\ell_\alpha):
u \in [\{0,\dotsc,t_m\}]^{t_m}\}
\endalign
$$ 
\medskip

$\quad$ [by dividing the sum $\dsize \sum_{\ell < m+1}$ into
$\dsize \sum_{\ell < m}$ and $\dsize \sum_{\ell =m}$ and noting what
$\nu^{n,k}_{\eta \char 94 \langle u \rangle}(m)$ is]
\bigskip

$= \dsize \sum 
\biggl\{ \bigl( \dsize \prod^{k-1}_{\ell =m+1}(t_\ell +1) \bigr)
\bigl( \dsize \sum_{\ell <m} \,\,
\dsize \sum_{\alpha \in \nu^{n,k}_{\eta \char 94 \langle u \rangle}(\ell)} 
p^\ell \, x^\ell_\alpha \bigr): u \in [\{0,\dotsc,t_m\}]^{t_m} \biggr\}$
\bigskip
$\qquad + \dsize \sum 
\biggl\{ \bigl( \dsize \prod^{k-1}_{\ell =m+1}(t_\ell +1) 
\bigr) \,\, \dsize \sum_{\alpha \in OP_{\omega,A_\eta}(u)} p^m \, x^m_\alpha: 
u \in [\{0,\dotsc,t_m\}]^{t_m} \biggr\} =$
\medskip

$\qquad \qquad \qquad \quad$[in the second sum, we collect
together the terms with $x^m_\alpha$]
\bigskip

$\qquad = \dsize \sum 
\biggl\{ \bigl( \dsize \prod^{k-1}_{\ell =m+1}(t_\ell +1) 
\bigr) \bigl( \dsize \sum_{\ell < m} \,\, 
\dsize \sum_{\alpha \in \nu^{n,k}_\eta(\ell)} p^\ell \, x^\ell_\alpha \bigr): 
u \in [\{0,\dotsc,t_m\}]^{t_m} \biggr\}$

$\qquad + \dsize \sum
\biggl\{ \bigl( \dsize \prod^{k-1}_{\ell =m+1}(t_\ell +1) 
\bigr)(p^m \, x^m_\alpha)|\{u:u \in [\{0,\dotsc,t_m\}]^{t_m} 
\text{ and } |\alpha \cap A_\eta| \in u\}|:$

$\qquad \qquad \qquad \alpha \text{ is a member of } A_\eta,
\text{ moreover } |\alpha \cap A_\eta| \le t_m \biggr\}$
\bigskip
$\qquad = \bigl( \dsize \prod^{k-1}_{\ell =m+1}(t_\ell +1) \bigr)
\bigl( \dsize \sum_{\ell <m} \,\,\dsize \sum_{\alpha \in 
\nu^{n,k}_\eta(\ell)} 
p^\ell \, x^\ell_\alpha \bigr) \times |\{u:u \in [\{0,\dotsc,t_m\}]^{t_m}\}|$ 
\bigskip
$\qquad + \dsize \sum
\biggl\{ \bigl( \dsize \prod^{k-1}_{\ell =m+1}(t_\ell +1) 
\bigr)(p^m \, x^m_\alpha) \cdot ((t_m +1)-1):\alpha \in A_\eta,|\alpha \cap
A_\eta| \le t_m \biggr\} =$  
\medskip

$\qquad \qquad \quad \qquad \qquad \qquad \qquad \qquad \quad$ 
[remember $t_m=p$ and $p^{m+1} \, x^m_\alpha=0$]
\bigskip

$\qquad \qquad =(t_m+1) \bigl( \dsize \prod^{k-1}_{\ell=m+1}(t_\ell +1) \bigr)
\dsize \sum_{\ell <m} \,\, \dsize \sum_{\alpha \in \nu^{n,k}_\eta(\ell)} 
p^\ell \, x^\ell_\alpha + 0$
\bigskip

$\qquad \qquad = \bigl( \dsize \prod^k_{\ell =m}(t_\ell +1) \bigr)
\bigl( \dsize \sum_{\ell <m} \,\, \dsize \sum_{\alpha \in 
\nu^{n,k}_\eta(\ell)}p^\ell \, x^\ell_\alpha \bigr)$.
\bigskip

Hence we have finished the proof of $\boxtimes$. \nl
So, as for $n \in \Gamma,B_n$ serves for every $k \in (n,\omega)$, 
if $u_1,u_2 \in [B_n]^{t_n}$ are distinct then 
$y_{\langle u_1 \rangle} - y_{\langle u_2 \rangle}$ contradicts $(*)$.
\hfill$\square_{\scite{2.4}}$ 
\enddemo
\bn
Recall
\definition{\stag{2.4A} Definition}  1) Let $I$ be an ideal on $\kappa$ (or
just $I \subseteq {\Cal P}(\kappa)$ closed downward, $I^+ = {\Cal P}(\kappa)
\backslash I$), then we let:

$$
\align
\bold U_I(\lambda) = \text{ Min} \bigl\{ |{\Cal P}|:&{\Cal P} \subseteq
[\lambda]^{\le \kappa} \text{ and for every } f \in {}^\kappa \lambda\\
  &\text{for some } a \in {\Cal P} \text{ we have } \{i < \kappa:f(i) \in a
\} \in I^+ \bigl\}.
\endalign
$$
\mn
2) For $\sigma \le \theta \le \mu \le \lambda$ let cov$(\lambda,\mu,\theta,
\sigma) = \text{ Min}\{\lambda + |{\Cal P}|:{\Cal P} \text{ is a family of
subsets of } \lambda$ each of cardinality $\mu \text{ such that any }
X \subseteq \lambda \text{ of cardinality } < \theta$ is included in the
union of $< \sigma$ members of ${\Cal P}\}$.
\enddefinition
\bigskip

\proclaim{\stag{2.5} Claim}  1) For every $\lambda \ge \beth_\omega$, for some
$\theta < \beth_\omega$ for every $\mu \in (\beth_{p-1}(\theta),\beth_\omega)$
we have (letting $\mu_n = \mu$) \, 
$\bold U_{I^1_{\bar \mu,\theta}}(\lambda) = \lambda$ (hence 
$\bold U_{I^0_{\bar \mu}}(\lambda) = \lambda$). \nl
2)  If cf$(\lambda) > \beth_\omega$, then for some $\theta < \beth_\omega$,
for every $\mu \in (\beth_{p-1}(\theta),\beth_\omega)$ and $\lambda' <
\lambda$ we have $\bold U_{I^1_{\bar \mu,\theta}}(\lambda') < \lambda$.
\endproclaim
\bigskip

\demo{Proof}  By \scite{2.4}, $I_{\bar \mu,\theta}$ is a $\theta$-complete
proper ideal on a set of cardinality $\mu^{\aleph_0}$, for any $\mu,\theta$
as in the assumptions.  By \cite{Sh:460} for each $\lambda' \le \lambda$ for 
some $\theta = \theta[\lambda'] < \beth_\omega$ for every
$\mu \in (\theta,\beth_\omega)$ we have cov$(\lambda',\mu^+,\mu^+,\theta) =
\lambda'$, i.e. there is ${\Cal P}_\mu \subseteq [\lambda']^\mu$ of 
cardinality $\le \lambda'$ such that: if
$Y \in [\lambda']^{\le \mu}$ then $Y$ is included in the union of $< \theta$
members of ${\Cal P}_\mu$.  As $I^1_{\mu,\theta}$ is a $\theta^+$-complete
ideal on a set of cardinality $\mu$ it follows that 
$\bold U_{I^1_{\mu,\theta}}(\lambda') \le \lambda' \times |{\Cal P}_\mu| = 
\lambda'$ (and trivially $\bold U_{I^1_{\bar \mu,\theta}}(\lambda) \ge 
\lambda$).  For some $\theta < \beth_\omega$, for arbitrarily large 
$\lambda' < \lambda,\theta[\lambda'] \le \theta$ and $\theta[\lambda] \le
\theta$; clearly we are done.  \hfill$\square_{\scite{2.5}}$
\enddemo
\bigskip

\demo{\stag{2.6} Conclusion}  If $\beth_\omega \le \mu^+ < \lambda =
\text{ cf}(\lambda) < \mu^{\aleph_0}$, \underbar{then} in 
${\frak K}^{\text{rs(p)}}_\lambda$ there is no universal member.
\enddemo
\bigskip

\demo{Proof}  By \scite{2.5} and \cite[5.9]{Sh:552}.
\enddemo
\bigskip

\noindent
Moreover
\proclaim{\stag{2.7} Claim}  Assume
\medskip
\roster
\item "{$(a)$}"  $\dsize \prod_{\ell < \omega} \kappa_\ell < \mu <
\lambda = {\text{\rm cf\/}}(\lambda) \le \lambda' < \mu^{\aleph_0}$
\smallskip
\noindent
\item "{$(b)$}"  $\mu^+ < \lambda$ or at least for some ${\Cal P}$
\sn
{\roster
\itemitem{ $(*)_{\Cal P}$ }  $|{\Cal P}| = \lambda \and
(\forall a \in {\Cal P})(a \subseteq \lambda \and {\text{\rm otp\/}}(a) = 
\mu)$ \nl

$\qquad \qquad \quad \and
(\forall E)(E \text{ a club of } \lambda \rightarrow (\exists a \in
{\Cal P})(a \subseteq E))$
\endroster}
\item "{$(c)$}"  $\lambda' = \bold U_{I^0_{\bar \kappa}}(\lambda)
< \mu^{\aleph_0}$ where $\bar \kappa = \langle \kappa_\ell:\ell < \omega
\rangle$ and note that $I^0_{\bar \kappa}$ depends on the prime $p$.
\endroster
\medskip

\noindent
\underbar{Then} we can find reduced separable abelian $p$-groups, 
$G_\alpha \in {\frak K}^{\text{rs(p)}}_\lambda$ for $\alpha < \mu^{\aleph_0}$ 
such that for every reduced separable abelian $p$-group $G$ of 
cardinality $\lambda'$ we have:
\mr
\item "{{}}"  some $G_\alpha$ is not embeddable into $G$; also the number
of ordinals $\alpha < \mu^{\aleph_0}$ such that ?
\ermn
Moreover, each $G_\alpha$ is slender, i.e. there is no homomorphism from
$\Bbb Z^\omega$ into $G_\alpha$ with range of infinite rank.
\endproclaim
\bigskip

\demo{Proof}  Same proof as that of \cite[5.9]{Sh:552},
\cite[7.5]{Sh:552}.
\enddemo
\newpage

\head {\S3 Non-existence of universals for $\aleph_1$-free abelian groups}
\endhead  \resetall
\bigskip

The first section dealt with ${\frak K}^{\text{rtf}}_\lambda$ improving
\cite{Sh:552}.  But the groups used there are ``almost divisible".  
So what occurs if we replace ${\frak K}^{\text{rtf}}$ by a variant 
avoiding this?  We suggest to
consider the $\aleph_1$-free abelian groups where type arguments like in
\S1 break down.  
So the proof of \cite{Sh:552} 
becomes relevant and it is natural to improve it as in \S2 (which deals with
${\frak K}^{\text{rs(p)}}$), for diversity we use a stronger ideal.  We have 
not looked at the problem for $\aleph_1$-free abelian groups 
of cardinality $\lambda$ when $\aleph_0 < \lambda < 2^{\aleph_0}$". \nl
So we concentrate here on torsion free (abelian) groups.
\bigskip

\definition{\stag{3.0} Definition}  1) Let $\bar t = \langle t_\ell:\ell < 
\omega \rangle,2 \le t_\ell < \omega$.  For abelian group $H$, the
$\bar t$-valuation is

$$
\|x\|_{\bar t} = \text{ Inf}\{2^{-m}:\dsize \prod_{\ell <m} t_\ell
\text{ divide } x \text{ (in } G)\}.
$$
\mn
This is a semi-norm.  Remember $d_{\bar t}(x,y) = \|x-y\|_{\bar t}$. 
This semi-norm induces a topology which we call the $\bar t$-adic topology.

If $t_\ell = p$ for $\ell < \omega$ we may write $p$ instead of $\bar t$.
\nl
2)  Let $c \ell_{\bar t}(A,H)$ be the closure of $A$ in $H$ under the 
$\bar t$-adic topology. \newline
Let $PC_H(X)$ be the pure closure of $X$ in $H$.  Moreover $PC^p_H(X)$ is the
$p$-adic closure in $H$ of the subgroup of $H$ which $X$ generates. \nl
3) Let ${\frak K}^{\text{rtf}}[\bar t]$ be the class of $\bar t$-reduced 
torsion free abelian groups, i.e. the $G \in {\frak K}^{\text{rtf}}$ such that
$\dbca_{n < \omega} \bigl(\dsize \prod_{i < n} t_i \bigr)G = \{0\}$ hence
$\|-\|_{\bar t}$ induces a Hausdorff topology. \nl
(Inversely if $G$ is torsion free with the $\bar t$-adic topology Hausdorff
then $G \in {\frak K}^{\text{rtf}}[\bar t]$.) \nl
4) If the $\bar t$-adic topology is Hausdorff, then $G^{[\bar t]}$ is
the completion of $G$ by $\|-\|_{\bar t}$. \nl
If $t_\ell = 2 + \ell$, this is the $\Bbb Z$-adic completion. 
\enddefinition
\bn
The following continues the analysis in \cite[1.1]{Sh:552} (which deals with
${\frak K}^{\text{rs(p)}}$) \nl
\cite[1.5]{Sh:552} (which deals with ${\frak K}^{\text{rtf}}$).
\bigskip

\definition{\stag{3.0A} Definition}  We say $G$ has $\bar t$-density $\mu$ 
if it has a pure subgroup of cardinality $\le \mu$ which is $\bar t$-dense, 
i.e. dense in the $\bar t$-adic topology, but has no such subgroup of 
cardinality $< \mu$.
\enddefinition
\bigskip

\proclaim{\stag{3.1} Proposition}  Suppose that
\medskip
\roster
\item "{$(\alpha)$}"  $\mu \le \lambda \le \mu^{\aleph_0}$
\smallskip
\noindent
\item "{$(\beta)$}"  $G$ is an $\aleph_1$-free abelian group, $|G| =
\lambda$
\sn
\item "{$(\gamma)$}"  $\bar t$ is as in \scite{3.0} such that
$(\forall \ell)(\exists m > \ell)$ ($t_\ell$ divides $t_m$).
\endroster
\medskip

\noindent
\ub{Then} there is an $\aleph_1$-free group $H$ such that 
$G \subseteq H,|H| = \lambda$ and $H$ has $\bar t$-density $\mu$.
\endproclaim
\bigskip

\demo{Proof}    Choose $\lambda_n < \mu$ for $n < \omega$ such that 
$\dsize \prod_{n < \omega} \lambda_n \ge \lambda,\mu \ge 
\dsize \sum_{n<\omega} \lambda_n,2\lambda_n < \lambda_{n+1}$ (so $\lambda_n 
> 0$ may be finite).  
Let $\{x_i:i < \lambda\}$ list the elements of $G$.  
Let $\lambda'_{n+1} = \lambda_{n+1},\lambda'_0 = \mu$ if 
$\mu > \dsize \sum_{n < \omega} \lambda_n$ and
$\lambda'_n = \lambda_n$ for all $n$ if 
$\mu = \dsize \sum_{n < \omega} \lambda_n$.
Let $\eta_i \in \dsize \prod_{n < \omega} \lambda'_n$ for $i < \lambda$
be distinct such that $\eta_i(n+1) \ge \lambda_n$ and $i \ne j \Rightarrow 
(\exists m)(\forall n)[m \le n \Rightarrow \eta_i(n) \ne \eta_j(n)]$, and
$\{\eta_i(0):i < \lambda\} = \lambda'_0$.
Without loss of generality $\mu = \{\eta_i(n):i < \lambda,n < \omega\}$.
Let $H$ be generated by $G,x^m_i \, (\text{for } i < \lambda'_m,m < \omega),
y^n_i \, (\text{for } i < \lambda,n < \omega)$ freely except for
\medskip
\roster
\item "{$(a)$}"  the equations of $G$
\smallskip
\noindent
\item "{$(b)$}"  $y^0_i = x_i \,(\in G)$
\smallskip
\noindent
\item "{$(c)$}"  $t_n y^{n+1}_i + y^n_i = x^n_{\eta_i(n)}$.
\endroster
\enddemo
\bigskip

\noindent
\underbar{Fact A}:  $H$ extends $G$ and is torsion free. 
\bigskip

\demo{Proof}   $H$ can be embedded into the divisible hull of 
$G \times F$, where $F$ is the abelian group generated freely by 
$\{x^m_\alpha:m < \omega \text{ and } \alpha < \lambda'_m\}$.
\enddemo\bigskip

\noindent
\underbar{Fact B}:  $H$ is $\aleph_1$-free and moreover $H/G$ is
$\aleph_1$-free.
\bigskip

\demo{Proof}  Let $K$ be a countable pure subgroup of $H$.  Now without loss
of generality $K$ is generated by
\medskip
\roster
\item "{$(i)$}"  $K_1 = \{x_i:i \in I \}$ is a pure subgroup of $G$, where 
$I$ is some countably infinite subset of $\lambda$, and so 
$G \supseteq K_1$,
\smallskip
\noindent
\item "{$(ii)$}"  $y^m_i,x^n_j$ for $i \in I,m < \omega$ and $(n,j) \in J$,
where $J \subseteq \omega \times \lambda$ is countable and

$$
i \in I,n < \omega \Rightarrow (n,\eta_i(n)) \in J.
$$
\endroster
\medskip

\noindent
Moreover, the equations holding among those elements are deducible from
the equations of the form
\medskip
\roster
\item "{$(a)^-$}"  equations of $K_1$
\smallskip
\noindent
\item "{$(b)^-$}"  $y^0_i = x_i$ for $i \in I$
\smallskip
\noindent
\item "{$(c)^-$}"  $t_n y^{n+1}_i + y^n_i = x^n_{\eta_i(n)}$ for 
$i \in I,n < \omega$.
\endroster
\medskip

\noindent
We can find $\langle k_i:i < \omega \rangle$ such that
$[i \ne j \and i \in I \and j \in I \and n \ge k_i \and n \ge k_j 
\and i \ne j \Rightarrow \eta_i(n) \ne \eta_j(n)]$.
\mn

Now we know that $K_1$ is free (being a countable subgroup of $G$), and it
suffices to prove that $K/K_1$ is free.  But $K/K_1$ is freely generated by
\nl
$\{y^n_i:i \in I \text{ and } n > k_i\} \cup \{x^n_\alpha:(n,\alpha) \in J
\text{ but for no } i \in I \text{ do we have}$ \nl
$n > k_i,\eta_i(n) = \alpha\}$.
So $K$ is free.
\enddemo
\bn
\ub{Fact C}:  $H_0 = \langle x^n_i:n < \omega,i < \lambda'_n \rangle_H$
satisfies:
\mr
\item "{$(a)$}"  $i < \lambda \Rightarrow d_{\bar t}(x_i,H_0) = \text{ inf}
\{d_{\bar t}(x_i,z):z \in H_0\}=0$
\sn
\item "{$(b)$}"  $x \in G \Rightarrow d_{\bar t}(x,H_0) = 0$
\sn
\item "{$(c)$}"  $x \in H \Rightarrow d_{\bar t}(x,H_0) = 0$.
\endroster
\bigskip

\demo{Proof}  First note that
\mr
\item "{$(*)_1$}"  $Y = \{x \in H:d_{\bar t}(x,H_0)=0\}$ \nl
\ermn
is a subgroup of $H$.  Also for every $i < \lambda$ and $n$
\mr
\item "{$(*)_2$}"  $y^n_i = x^n_{\eta_i(n)} + t_n y^{n+1}_i = 
x^n_{\eta_i(n)} + t_n x^{n+1}_{\eta_i(n+1)} + 
t_n t_{n+1} y^{n+2}_i$ \nl

$= \dsize \sum^m_{k=n} \bigl( \dsize \prod^{k-1}_{\ell =n} 
t_\ell \bigr) x^k_{\eta_i(k)} 
+ \bigl( \dsize \prod^k_{\ell =n} t_\ell \bigr) y^{k+1}_i$
\ermn
(prove by induction on $m \ge n$), and note that as $(\forall \ell)(\exists
m > \ell)(t_i$ divides $t_m$) necessarily $(\forall \ell)(\exists^\infty m)
(t_\ell$ divides $t_m$) hence $(\forall k)(\exists^\infty m)
(\dsize \prod_{i \le \ell} t_\ell$ divides $\dsize \prod^m_{i=k} t_i$).  Now 
$(*)_2$ implies
\mr
\item "{$(*)_3$}"  $y^n_i \in Y$.
\ermn
But $x_i = y^0_i$ and hence clause (a) holds, so as $\{x_i:i < \lambda\}$ is
dense in $G$ also clause (b) holds.  So $G \subseteq Y$ (by clause (b)), and
$x^n_\alpha \in Y$ (as $H_0 \subseteq Y$ and the choice of $H_0$) and
$y^n_i \in Y$ (by $(*)_3$). \nl
By $(*)_1$ clearly $Y=H$, as required in clause (c).
\enddemo
\bn
\ub{Fact D}:  $|H|=\lambda$. 
\bn
\ub{Fact E}:  The $\bar t$-density of $H$ is $\mu$. 
\bigskip

\demo{Proof}  It is $\le \mu$ as $H_0$ has cardinality $\mu$ and 
is $\bar t$-dense in $H$, it is $\ge \mu$, as we now show.  
Otherwise let $X^* \subseteq H$ be $\bar t$-dense with $|X^*| < \mu$, so 
for some $n,|X^*| < \lambda'_n$.

Define a function $h$ with domain the generators of $H$ listed above, into
$H$.  Let $h(x) = 0 \text{ if } x \in G;h(x^m_\alpha) = 0 
\text{ if } m \ne n,\alpha < \lambda_m;h(x^m_\alpha) = 
x^m_\alpha \text{ if } m = n,\alpha < \lambda_m;h(y^m_i) = 0 \text{ if } 
m > n;h(y^m_i) = -x^n_{\eta_i(n)} \text{ if } m = n;h(y^m_i) = 
\dsize \prod^{n-1}_{\ell = m} t_\ell \times h(y^n_i) \text{ if } m < n$.
\nl
This function preserves the equations defining $H$ and hence induces a
homomorphism $\hat h$ from $H$ onto $\langle \text{Rang}(h) \rangle_H =
\langle \{x^n_\alpha:\alpha < \lambda'_n\} \rangle_H$.
Clearly $h(h(x)) = h(x)$ for the generators hence $\hat h \circ \hat h = 
\hat h$.  Hence $\langle \{x^n_\alpha:\alpha < \lambda'_n\} \rangle_H$ 
is a direct summand
of $H$ and hence the $d_{\bar t}$-density of $H$ is at least the
$d_{\bar t}$-density of $\langle \{x^n_\alpha:\alpha < \lambda'_n\} 
\rangle_H$ which is $\lambda'_n$.  \hfill$\square_{\scite{3.1}}$
\enddemo
\bn
We define variants of Definition \scite{2.1}.
\definition{\stag{3.2} Definition}  For $\bar \lambda = \langle \lambda_\ell:
\ell < \omega \rangle,\bar t = \langle t_\ell:\ell < \omega \rangle,
2 \le t_\ell < \omega$, we let

$$
\align
J^5_{\bar t,\bar \lambda} = \biggl\{ X \subseteq \dsize \prod_{\ell < \omega}
[\lambda_\ell]^{t_\ell}:&\text{ we cannot find } m(*) < \omega,
\bar Y = \langle Y_m:m < \omega \text{ and } m \ge m(*) \rangle,\\
  &\bar A^m = \langle A_\eta:\eta \in Y_m \rangle \text{ such that:} \\
  &\,(a) \qquad Y_m \subseteq \dsize \prod_{\ell < m} [\lambda_\ell]
^{t_\ell} \\
  &\,(b) \qquad Y_{m(*)} \subseteq \dsize \prod_{\ell < m(*)} [\lambda_\ell]
^{t_\ell} \text{ is a singleton} \\
  &\,(c) \qquad \langle A_\eta:\eta \in Y_m \rangle \text{ is a sequence of
pairwise disjoint} \\
  &\qquad \quad\text{ subsets of } \lambda_m \text{ of order type } \omega \\
  &\,(d) \qquad Y_{m+1} = \{\eta \char 94 \langle u \rangle:\eta \in Y_m
\text{ and } u \in [A_\eta]^{t_m}\} \\
  &\,(e) \qquad Y_m \subseteq \{\nu \restriction m:\nu \in X\} \biggr\},
\endalign
$$
\bn
$J^6_{\bar t,\bar \lambda} \text{ is defined similarly but } m(*)=0$,
\bigskip

$$
J^\ell_{\bar t,\bar \lambda,< \theta} = \biggl\{X:\text{ for some } \alpha <
\theta \text{ and } X_\beta \in J^\ell_{\bar t,\bar \lambda} \text{ for }
\beta < \alpha \text{ we have } X \subseteq \dsize \bigcup_{\beta < \alpha} 
X_\beta \biggr\}.
$$
Also let $J^\ell_{\bar t,\bar \lambda,\theta} = 
J^\ell_{\bar t,\bar \lambda,< \theta^+}\,\,$.
\enddefinition 
\bigskip

\proclaim{\stag{3.3} Claim}  1)  $J^{i(1)}_{\bar t,\bar \lambda,< \theta_1}
\subseteq J^{i(2)}_{\bar t,\bar \lambda,< \theta_2}$ when $\theta_1 \le
\theta_2,i(1) \le i(2)$ are mong 4,5,6. \newline
2) $J^i_{\bar t,\bar \lambda,\theta}$ is a $\theta^+$-complete ideal for
$i = 4,5,6$. \newline
3) If $\lambda_\ell \ge \beth_{t_\ell-1}(\theta)$ then the ideal 
$J^i_{\bar t,\bar \lambda,\theta}$ is proper for $i =4,5,6$.
\endproclaim
\bigskip

\demo{Proof}  1), 2) Easy. \newline
3) As in \scite{2.4}.  \hfill$\square_{\scite{3.3}}$
\enddemo
\bigskip

\definition{\stag{3.4} Definition}  For $\bar \lambda = \langle \lambda_\ell:
\ell < \omega \rangle,\bar t = \langle t_\ell:\ell < \omega \rangle$ 
such that $2 \le t_\ell < \omega$ we define
\medskip
\roster
\item "{$(A)$}"  $B^{\text{rtf}}_{\bar t,\bar \lambda}$ is the 
free (abelian) group
generated by $\{x^m_\alpha:m < \omega,\alpha < \lambda_m\}$. 
\smallskip
\noindent
\item "{$(B)$}"  Let $B^{\text{rtf}}_{\bar t,\bar \lambda,n}$ 
be the subgroup of $B^{\text{rtf}}_{\bar t,\bar \lambda}$ generated by 
$\{x^m_\alpha:m < n$ and $\alpha < \lambda_m\}$
\smallskip
\noindent
\item "{$(C)$}"  $G^{\text{rtf}}_{\bar t,\bar \lambda}$ is the pure closure in
$(B^{\text{rtf}}_{\bar t,\bar \lambda})^{[\bar t]}$ of the subgroup of 
$(B^{\text{rtf}}_{\bar t,\bar \lambda})^{[\bar t]}$ generated by
$$
B^{\text{rtf}}_{\bar t,\bar \lambda} \cup \bigl\{ \dsize \sum_{m < \omega}
\bigl( \dsize \prod_{\ell < m} t_\ell \bigr) 
(x^\ell_{(\eta(\ell))(1)} - x^\ell_{(\eta(\ell))(0)}):\eta \in
\dsize \prod_{\ell < \omega}[\lambda_\ell]^2 \bigr\}
$$
(here we use that if e.g. 
$\eta(\ell) = \{\alpha,\beta\},\alpha < \beta$ then 
$(\eta(\ell))(1) = \beta,(\eta(\ell))(0) = \alpha$.
\smallskip
\noindent
\item "{$(D)$}"  Let $\bar B^{\text{rtf}}_{\bar t,\bar \lambda} = \langle
B^{\text{rtf}}_{\bar t,\bar \lambda,n}:n < \omega \rangle$.
\endroster
\enddefinition
\bigskip

\definition{\stag{3.6} Definition}  Assume
\medskip
\roster
\item "{$\boxtimes^{\bar t}_{H,\bar H}$}"  $\bar H = \langle H_n:
n < \omega \rangle$ is an increasing sequence of abelian subgroups of $H$,
such that $\dbcu_{n < \omega} H_n$ is dense in $H$ by the $\bar t$-adic
topology.
\endroster
\medskip

\noindent
Then we let

$$
\align
I^{4,\bar t}_{H,\bar H} = \bigl\{X \subseteq H:&\text{ for some } n < \omega,
\text{ the intersection of the } \bar t \text{-adic closure of } 
PC(X)_H \text{ in } H,\\
  &c \ell_{\bar t}(PC(X)_H,H)\text{ with } 
\dsize \bigcup_{\ell < \omega} H_\ell \text{ is a subset of } H_n \bigr\}
\endalign
$$

$$
I^{4,\bar t}_{H,\bar H,<\theta} = \bigl\{X \subseteq H:\text{ for some } 
\alpha < \theta \text{ and } X_\beta \in I^{4,\bar t}_{H,\bar H} 
 \text{ for } \beta < \alpha \text{ we have } X \subseteq \dsize
\bigcup_{\beta < \alpha} X_\beta \bigr\}
$$

$\qquad \qquad I^{4,\bar t}_{H,\bar H,\theta} = 
I^{\bar t}_{H,\bar H,< \theta^+}$.
\enddefinition
\bigskip

\definition{\stag{3.7} Definition}  Assume $\bar t = \langle t_\ell:\ell < 
\omega \rangle,2 \le t_\ell < \omega$, and 
\medskip
\roster
\item "{$\boxtimes^{\bar t}_{H,\bar H}$}"  $H$ is Hausdorff in the $\bar t
\restriction [k,\omega)$-topology for each $k < \omega$ where $t \restriction
[k,\omega) = \langle t_{k + \ell}:\ell < \omega \rangle$.  Further
$\bar H = \langle H_n:n < \omega \rangle$ is an increasing sequence of 
abelian groups, $\dbcu_{n < \omega} H_n \subseteq H$ is dense in the
$\bar t \restriction [k,\omega)$-adic topology for each $k < \omega$.
\endroster
\medskip

\noindent
Then we let

$$
\align
I^{5,\bar t}_{H,\bar H} = \bigl\{X \subseteq H:&\text{ for some } n(*) < 
\omega, \text{ for every } n \in (n(*),\omega) \text{ there is no}
\tag "{$1)$}" \\
  &\,y \in H_{n+1} \text{ such that: } d_{\bar t \restriction [n,\omega)}
(y,PC(\langle X \rangle)) = 0 \text{ but } 
d_{\bar t \restriction [n,\omega)}(y,H_n) > 0 \bigr\}
\endalign
$$

$$
I^{5,\bar t}_{H,\bar H,<\theta} = \bigl\{X:\text{ there are } \alpha < 
\theta \text{ and } X_\beta \in I^{5,\bar t}_{H,\bar H} \text{ for }
\beta < \alpha \text{ such that } X \subseteq \dbcu_{\beta < \alpha} 
X_\beta \bigr\}.
$$

Moreover $I^{5,\bar t}_{H,\bar H,\theta} = I^{5,t}_{H,\bar H,< \theta^+}$.
\mn
2) $I^{6,\bar t}_{H,\bar H}$ (and $I^{6,\bar t}_{H,\bar H,<\theta},
I^{6,\bar t}_{H,\bar H,\theta}$) are defined similarly except that we demand
$n(*)=0$. \nl
\sn
3) $I^{i,\text{rtf}}_{\bar t,\bar \lambda}$ means
$I^{i,\bar t}_{G^{\text{rtf}}_{\bar t,\bar \lambda},
\bar B^{\text{rtf}}_{\bar t,\bar \lambda}}$ where 
$\bar B^{\text{rtf}}_{\bar t,\bar \lambda} =
\langle B^{\text{rtf}}_{\bar t,\bar \lambda,n}:n < \omega \rangle$. 
\enddefinition

\proclaim{\stag{3.8} Claim}  For $\bar \lambda,\bar t$ as in \scite{3.2}
\medskip
\roster
\item "{$(a)$}"  we have $\boxtimes^{\bar t}_{G^{\text{rtf}}
_{\bar t,\bar \lambda},\bar B^{\text{rtf}}_{\bar t,\bar \lambda}}$ 
(from \scite{3.7})
\smallskip
\noindent
\item "{$(b)$}"  $G^{\text{rtf}}_{\bar t,\bar \lambda}$ is 
$\aleph_1$-free; moreover $G^{\text{rtf}}_{\bar t,\bar \lambda}/
B^{\text{rtf}}_{\bar t,\bar \lambda,n}$ is $\aleph_1$-free for each 
$n < \omega$
\smallskip
\noindent
\item "{$(c)$}"  $I^{i,\text{rtf}}_{\bar t,\bar \lambda,\theta}$ are
$\theta^+$-complete ideals for $i = 4,5,6$
\sn
\item "{$(d)$}"  if $\boxtimes^{\bar t}_{H,\bar H}$ (from \scite{3.7}) and 
$i \in \{4,5,6\}$ then $I^{i,\bar t}_{H,\bar H,\theta}$ is a 
$\theta^+$-complete ideal.
\endroster
\endproclaim
\bigskip

\demo{Proof}  Straightforward (for (6), use an argument similar to that of
\scite{3.1}).
\enddemo
\bigskip

\noindent
The following lemma connects the combinatorial ideals defined above
and the more algebraic ideals defined in \scite{3.6}.
\proclaim{\stag{3.9} Claim}  1) Assume
\medskip
\roster
\item "{$\boxtimes_1$}"  $\bar t = \langle t_\ell:\ell < \omega \rangle,
2 \le t_\ell < \omega$
\smallskip
\noindent
\item "{$\boxtimes_2$}"  $\bar \lambda = \langle \lambda_\ell:\ell <
\omega \rangle$, and $\lambda_\ell > \beth_1(\theta)$ for $\ell <
\omega$.
\endroster
\medskip

\noindent
\ub{Then} the ideal $I^{i,\text{rtf}}_{\bar t,\bar \lambda,\theta}$ is 
proper for $i = 4,5,6$. \nl
2) Assume $\boxtimes_1$ and 
\mr
\item "{$\boxtimes'_2$}"  $\bar \lambda = \langle \lambda_\ell:
\ell < \omega \rangle,\lambda_\ell = \aleph_0,\theta = \aleph_0$. 
\ermn
\ub{Then} the ideal $I^{i,\text{rtf}}_{\bar t,\bar \lambda,\theta}$ is proper.
\endproclaim
\bigskip

\demo{Proof}  1) If not, we can find 
$X_\alpha \subseteq L =: G^{\text{rtf}}_{\bar t,
\bar \lambda}$ for $\alpha < \theta$ such that 
$G^{\text{rtf}}_{\bar t,\bar \lambda}
= \dsize \bigcup_{\alpha < \theta} X_\alpha$ and $X_\alpha \in
I^{i,\text{rtf}}_{\bar t,\bar \lambda}$.  For $\alpha \le \omega$ and 
$\eta \in \dsize \prod_{\ell < \alpha} [\lambda_\ell]^2$ we let \nl 
$x_\eta = \dsize \sum_{m < \alpha} \bigl( \dsize \prod_{\ell < m} t_\ell
\bigr)(x^m_{(\eta(n))(1)} - x^m_{(\eta(n))(0)})$.  
\mn
As in the proof of \scite{2.4}, we apply \cite{RuSh:117}, \relax
\cite[Ch.XI, 3.5]{Sh:f} for the ideal $J_\ell = ERI^2_\theta(\lambda_\ell)$ (this ideal is, 
of course, $\theta^+$-complete and non-trivial as $\lambda_\ell > 2^\theta$).
\medskip

So we can find $\langle Y_m:m < \omega \rangle,\langle A_\eta:\eta \in
Y_m \rangle$ and $\alpha(*) < \theta$ such that
\medskip
\roster
\item "{$(a)$}"  $Y_m \subseteq \dsize \prod_{\ell < m}[\lambda_\ell]
^{t_\ell}$
\smallskip
\noindent
\item "{$(b)$}"  $Y_0$ is a singleton
\smallskip
\noindent
\item "{$(c)$}"  $A_\eta \in (J_{\ell g(\eta)})^+$ for $\eta \in Y_m$ (so
$A_\eta \subseteq [\lambda_{\ell g(\eta)}]^{t_{\ell g(\eta)}}$)
\smallskip
\noindent
\item "{$(d)$}"  $Y_{m+1} = \{\eta \char 94 \langle u \rangle:u \in A_\eta,
\eta \in Y_m\}$
\smallskip
\noindent
\item "{$(e)$}"  if $\eta \in Y_m$ then $\eta \in \{\nu \restriction m:
x_\nu \in X_{\alpha(*)}\}$.
\endroster
\medskip

\noindent
Note that we can demand that for all $\alpha < \beta$ from $A_\eta,
d(x^{\ell g(\eta)}_\beta - x^{\ell g(\eta)}_\alpha,\langle X_{\alpha(*)}
\rangle)$ is the same.
We now prove by induction on $k < \omega$ that
\mr
\item "{$(*)_k$}"  for any $m < \omega$, if $\eta \in Y_m$ and 
$A \subseteq A_\eta$ is infinite then for some infinite $A' \subseteq A$ for
any $\alpha < \beta$ from $A'$ we have \nl
$(\dsize \prod_{\ell < m} t_\ell)(x^m_\beta - x^m_\alpha) \in c \ell_{\bar t}
(\langle X_{\alpha(*)} \rangle,L) + 
(\dsize \prod_{\ell < m+k} t_\ell) L$. 
\ermn
For $k=0$ this is trivial: the element $(\dsize \prod_{\ell < m}
t_\ell)(x^m_\beta - x^m_\alpha)$ belongs to $(\dsize \prod_{\ell < m+k}
t_\ell)L$. \nl
For $k+1$, by $(*)_k$ as we can replace every $A_\eta$ by an infinite subset
without loss of generality for $m < \omega,A_\eta$ can serve as $A'$ in 
$(*)_k$.  Let $\eta \in Y_m$, and let $\gamma_0 < \gamma_1 < \gamma_2
< \ldots$ be in $A_\eta$.  So for each $j < \omega$, let
$\eta_j \in Y_{m+k+1}$ be such that $\eta_j \restriction m = \eta,\eta_j(m)
= \{\gamma_j,\gamma_{j+1}\}$.  By clause (e) above we know that there are 
$\nu_j$ such that 
$\eta_j \triangleleft \nu_j \in \dsize \prod_{\ell < \omega}[\lambda_\ell]^2$ 
and
\mr
\widestnumber\item{$(iii)$}
\item "{$(i)$}"  $x_{\nu_j} \in X_{\alpha(*)}$. \nl
Now by the definitions of $x_{\eta_j},x_{\nu_j}$
\sn
\item "{$(ii)$}"  $x_{\eta_j} = x_{\nu_j}$ mod$(\dsize \prod_{\ell < m+k+1}
t_\ell)L$ 
\sn
\item "{$(iii)$}"  if $\ell \in [m+1,m+k+1)$ and $j < \omega$ \ub{then} \nl
$x_{\eta_j \restriction (\ell + 1)} - x_{\eta_j \restriction \ell} \in
c \ell_{\bar t}(\langle X_{\alpha(*)} \rangle,L) + 
(\dsize \prod_{i < \ell + k} t_i) L \subseteq c \ell_{\bar t}
(\langle X_{\alpha(*)} \rangle,L) + (\dsize \prod_{i<m+k+1} t_i)L$ \nl
[why? by the induction hypothesis the difference is \nl
$(\dsize \prod_{i <m+\ell}t_i)(x^\ell_{(\eta_j(\ell))(1)} - 
x^\ell_{(\eta_j(\ell))(0)})$]
\sn
\item "{$(iv)$}"  $x_{\eta_j} - x_{\eta_j \restriction (m+1)} \in
c \ell_{\bar t}(\langle X_{\alpha(*)} \rangle,L) + 
(\dsize \prod_{i < m+k+1} t_i)L$ \nl
[why?  use (iii) for $\ell = m+1,\dotsc,m+k$, noting that $\ell g(\eta_j) =
m+k+1$.]
\sn
\item "{$(v)$}"  $x_{\eta_j \restriction (m+1)} \in c \ell_{\bar t}
(\langle X_{\alpha(*)} \rangle,L) + (\dsize \prod_{i < m+k+1} t_i)L$ \nl
[why?  by (i) + (ii) + (iv)]
\sn
\item "{$(vi)$}"  $\dsize \sum \bigl\{ x_{\eta_j \restriction (m+1)}:j <
\dsize \prod_{i < m+k+1} t_i \bigr\} \in c \ell_{\bar t}
(\langle X_{\alpha(*)} \rangle,L) + (\dsize \prod_{i < m+k+1} t_i)L$ \nl
[why?  by (v)]
\sn
\item "{$(vii)$}"  $x^m_{\gamma_{j(*)}} - x^m_{\gamma_0} \in c \ell_{\bar t}
(\langle X_{\alpha(*)} \rangle,L) + (\dsize \prod_{i < m+k+1} t_i))L$ 
for $j(*) = \dsize \prod_{i < m+k+1} t_i$ \nl
[why?  by (vi) because
\endroster
\bn
$\dsize \sum 
\bigl\{ x_{\eta_j \restriction (m+1)}:j < \dsize \prod_{i < m+k+1} t_i
\bigl\}$ 
\mn
$\qquad = \dsize \sum \bigl\{ x_{\eta_j \restriction m} + 
(\dsize \prod_{i < m+1} t_i)(x^m_{\gamma_{j+1}} - x^m_{\gamma_j}):j <
\dsize \prod_{i<m+k+1} t_i \bigr\}$
\mn
$\qquad = \dsize \sum \bigl\{ x_{\eta_j \restriction m}:j < 
\dsize \prod_{i < m+k+1} t_i \bigl\} + (\dsize \prod_{i<m+1} t_i)
\dsize \sum \bigl\{ (x^m_{\gamma_{j+1}} - x^m_{\gamma_j}):j <
\dsize \prod_{i<m+k+1} t_i \bigr\}$
\medskip

$\qquad \qquad \qquad \quad$ 
[as $\eta_j \restriction m$ does not depend on $j$ and
obvious arithmetic]
\mn
$\qquad = (\dsize \prod_{i < m+k+1} t_i) \cdot x_{\eta_{j(*)}} +
(\dsize \prod_{i<m+1} t_i)(x^m_{\gamma_{j(*)}} - x^m_{\gamma_0}) \in$
\mn
$\qquad \qquad \qquad (\dsize \prod_{i<m+1} t_i)(x^m_{\gamma_{j(*)}} -
x^m_{\gamma_0}) + (\dsize \prod_{i<m+k+1} t_i)L$]
\bn
\mr
\item "{$(viii)$}"  if $\rho \in Y_m$ and $\alpha < \beta$ are in $A_\eta$
then \nl
$(\dsize \prod_{i < m+1} t_i)(x^m_\beta - x^m_\alpha) \in c \ell_{\bar t}
(\langle X_{\alpha(*)} \rangle,L) + (\dsize \prod_{i<m+k+1} t_i)L$ \nl
[why?  by (vii) and the choice of the $Y_m,A_\eta(\eta \in Y_m,m < \omega)$.]
\ermn
So we have carried the induction on $k$.  \nl
2) Easier. \enddemo \hfill$\square_{\scite{3.9}}$
\bigskip

\proclaim{\stag{3.10} Claim}  Assume
\medskip
\roster
\item "{$\boxtimes_1$}"  $\bar t = \langle t_\ell:\ell < \omega \rangle$
and $2 \le t_\ell < \omega$
\smallskip
\noindent
\item "{$\boxtimes_2$}"  $\lambda_\ell > \beth_1(\theta)$
\smallskip
\noindent
\item "{$\boxtimes_3$}"  cov$(\lambda,\bigl( \dsize \prod_{\ell < \omega}
\lambda_\ell \bigr)^+, \bigl( \dsize \prod_{\ell < \omega} \lambda_\ell
\bigr)^+,\theta^+) \le \lambda$.
\endroster
\medskip

\noindent
\ub{Then} $\bold U_{J^6_{\bar t,\bar \lambda,\theta}}(\lambda) = \lambda$ and
$\bold U_{I^6_{\bar t,\bar \lambda,\theta}}(\lambda) = \lambda$.
\endproclaim
\bigskip

\demo{\stag{3.11} Conclusion}  For every $\lambda \ge \beth_\omega$
for some $\theta < \beth_\omega$, for every $\kappa \in (\beth_1(\theta),
\beth_\omega)$ for every $\lambda_n \in [\beth_\omega(\theta),
\kappa]$ we have

$$
\bold U_{I^6_{\bar t,\bar \lambda,\theta}}(\lambda) = \lambda = 
\bold U_{J^6_{\bar t,\bar \lambda,\theta}}(\lambda).
$$
\enddemo
\bigskip

\demo{Proof}  By the previous claim and \cite{Sh:460} (similar to
\scite{2.5}).  \hfill$\square_{\scite{3.11}}$
\enddemo
\bigskip

\proclaim{\stag{3.12} Claim}  Assume
\medskip
\roster
\item "{$(a)$}"  $\dsize \prod_{\ell < \omega} \lambda_\ell < \mu < \lambda
= \text{ cf}(\lambda) \le \lambda' \le \lambda'' < \mu^{\aleph_0}$
\smallskip
\noindent
\item "{$(b)$}"  $\mu^+ < \lambda$ or at least for some ${\Cal P}$,
{\roster
\itemitem{ $(*)_{\Cal P}$ }  $|{\Cal P}| = \lambda$ and
$(\forall a \in {\Cal P})(a \subseteq \lambda \and \text{{\rm otp\/}}
(a) = \mu)$ \nl
and $(\forall E)(E \text{ a club of } \lambda \rightarrow (\exists a \in
{\Cal P})(a \subseteq E))$
\endroster}
\item "{$(c)$}"  $\lambda'' = \bold U_{I^6_{\bar t,\bar \lambda}}(\lambda')
< \mu^{\aleph_0}$ where $t_m = \dsize \prod_{\ell < m} \ell!$ or at
least $\lambda'' = \bold U_{J^1_{\bar \lambda}}(\lambda')$
\smallskip
\noindent
\item "{$(d)$}"  cov$(\lambda'',\lambda^+,\lambda^+,\lambda) < \mu^{\aleph_0}$
or at least $\bold U_{\text{id}^a({\Cal P})}(\lambda'') < \mu^{\aleph_0}$
where ${\Cal P}$ satisfies $(*)_{\Cal P}$.
\endroster
\medskip

\noindent
\underbar{Then} we can find $\aleph_1$-free abelian groups $G_\alpha$ of
cardinality $\lambda$ for $\alpha < \mu^{\aleph_0}$ such that for every 
$\aleph_1$-free abelian group $G$ of cardinality $\lambda$ or just 
$G \in {\frak K}^{\text{rtf}}_\lambda$ we have:
\mr
\item "{{}}"  some $G_\alpha$ is not embeddable into $G$; also the number of
ordinals $\alpha < \mu^{\aleph_0}$ for which $G_\alpha$ is embeddable into 
$G$ is at most cov$(\lambda'',\lambda^+,\lambda^+,\lambda)$ (or 
$\le \bold U_{\text{id}^a({\Cal P})}(\lambda'')$ at least)
\endroster
\endproclaim
\bigskip

\demo{Proof}  Like \scite{2.7}, 
not that ``$\aleph_1$-free" implies $\|-\|_{\bar t}$ is a norm.
\enddemo
\bigskip

\demo{\stag{3.13} Conclusion}  If $\beth_\omega \le \mu^+ < \lambda =
\text{ cf}(\lambda) < \lambda^{\aleph_0}$ \underbar{then} in
${\frak K}^{\text{rtf}}_\lambda$ there is no member universal even just for
${\frak K}^{\aleph_1\text{-free}}_\lambda$.
\enddemo
\bn

\demo{Proof}  Straightforward.
\enddemo
\bigskip

\remark{\stag{3.14} Concluding Remarks}  1) We can in \scite{3.10} -
\scite{3.13} use \scite{3.9}(2) instead of \scite{3.9}(1). \nl
2) In \S2 we can use the parallel of \scite{3.9}.  Also we can deal with
the class of $R$-modules or some natural subclasses of it, we hope to return
this elsewhere.
\endremark
\bigskip

\remark{\stag{3.15} Remark}  If $\lambda = \aleph_0$ there is no universal
member in ${\frak K}^{\text{rtf}}_\lambda$. 
In fact for any $\bold Q \subseteq \bold P^*$ let $G_{\bold Q}$ be the 
subgroup of $\Bbb Q x \oplus {\underset p {}\to \bigoplus}
\{\Bbb Q x_p:p \in \bold P^* \backslash \bold Q\}$ generated by

$$
\align
\{p^{-n} x:p \in \bold Q\} &\cup \{q^{-n} x_p:p \in \bold P^* \backslash
\bold Q \text{ and } n < \omega, \text{ and } q \in \bold P^* \backslash 
\{p\}\} \\
  &\cup \{ p^{-n}(x-x_p):n < \omega \text{ and } p \in \bold Q\}.
\endalign
$$
\mn
So $G_{\bold Q} \in {\frak K}^{\text{rtf}}_{\aleph_0}$, and if $h$ embeds
$G_{\bold Q}$ into $G \in K^{\text{trf}}$ then $\bold P(h(x),G) = \bold Q$.
\endremark

\newpage
    
REFERENCES.  
\bibliographystyle{lit-plain}
\bibliography{lista,listb,listx,listf,liste}

\shlhetal
\enddocument

\bye